\documentclass[11pt, reqno]{amsart}
\usepackage{amsmath}

\usepackage{amsthm}
\usepackage{amssymb}
\usepackage{amsmath,amssymb,tikz-cd}
\usepackage{mathtools}
\usepackage{bbm}
\usepackage{comment}
\usepackage{tikz}

\usepackage{lipsum}

\usepackage[english]{babel}

\usepackage{booktabs}

\usepackage{bm}

\usepackage{hyperref}
\hypersetup{
   colorlinks = true,
   linkcolor = {red},
   citecolor = {blue},
}
\usepackage{cleveref}

\usepackage{gensymb}

\usepackage{eucal}

\usepackage{graphicx}
\graphicspath{ {./images/} }

\usepackage{geometry}
\usepackage[above,below,verbose,section]{placeins}

\usepackage[sc]{mathpazo}
\usepackage[T1]{fontenc}
\usepackage{microtype}
\usepackage{boxedminipage}

\usepackage{cite}

\newcommand{\R}{\mathbb{R}}

\newcommand{\N}{\mathbb{N}}
\newcommand{\eps}{\epsilon}

\newcommand{\KR}{\textrm{KR}}
\newcommand{\KL}{\textrm{KL}}
\newcommand{\dkl}{\mathcal{D}_{\textrm{KL}}}
\renewcommand{\varepsilon}{\epsilon}
\DeclareMathOperator*{\argmax}{arg\,max}
\DeclareMathOperator*{\argmin}{arg\,min}

\newtheoremstyle{note}%
{10pt}%
{3pt}%
{\itshape}%
{0pt}%
{\bfseries}%
{.}%
{0.5em}%
{}%

\theoremstyle{note}
\newtheorem{theorem}{Theorem}[section]
\newtheorem{corollary}[theorem]{Corollary}
\newtheorem{proposition}[theorem]{Proposition}
\newtheorem{lemma}[theorem]{Lemma}
\newtheorem{rem}[theorem]{Remark}

\newtheorem{assumption}{Assumption}

\newtheorem{definition}[theorem]{Definition}

\crefname{lemma}{Lemma}{Lemmas}
\crefname{theorem}{Theorem}{Theorems}
  \crefname{remark}{Remark}{Remarks}
  \crefname{proposition}{Proposition}{Propositions}
  \crefname{section}{Section}{Sections}
  \crefname{subsection}{Subsection}{Subsections}
  \crefname{equation}{}{}
  \Crefname{equation}{Equation}{Equations}
  \crefname{figure}{Figure}{Figures}
  \crefname{appendix}{Appendix}{Appendices}
\crefname{assumption}{Assumption}{Assumptions}
\crefname{enumi}{}{}

\usepackage{enumitem}
\makeatletter
\def\namedlabel#1#2{\begingroup
    #2%
    \def\@currentlabel{#2}%
    \phantomsection\label{#1}\endgroup
}
\makeatother

\usepackage{todonotes}

\makeatletter
\renewcommand{\paragraph}{%
  \@startsection{paragraph}{4}%
  {\z@}{1.0ex \@plus 1ex \@minus .2ex}{-1em}%
  {\normalfont\normalsize\bfseries}%
}
\makeatother

\title{Knothe-Rosenblatt maps via soft-constrained optimal transport}%

\author{\textsc{Ricardo Baptista}, \textsc{Franca Hoffmann}, \textsc{Minh Van Hoang Nguyen}, \textsc{Benjamin Zhang}}

\address{Department of Statistical Sciences and Vector Institute\\
University of Toronto\\
700 University Avenue\\
Toronto, ON M5G 1X6
}
\email{r.baptista@utoronto.ca}

\address{Computing + Mathematical Sciences\\
California Institute of Technology\\
1200 East California Boulevard\\
Pasadena, CA 91125}
\email{franca.hoffmann@caltech.edu}

\address{Computing + Mathematical Sciences\\
California Institute of Technology\\
1200 East California Boulevard\\
Pasadena, CA 91125}
\email{mvnguyen@caltech.edu}

\address{Division of Applied Mathematics\\
Brown University\\
182 George Street\\
Providence, RI 02912}
\email{benjamin\_zhang@brown.edu}



\begin{document}

\begin{abstract}
In the theory of optimal transport, the Knothe-Rosenblatt (KR) rearrangement provides an explicit construction to map between two probability measures by building one-dimensional transformations from the marginal conditionals of one measure to the other. The KR map has shown to be useful in different realms of mathematics and statistics, from proving functional inequalities to designing methodologies for sampling conditional distributions. 
It is known that the KR rearrangement can be obtained as the limit of a sequence of optimal transport maps with a weighted quadratic cost. We extend these results in this work by showing that one can obtain the KR map as a limit of maps that solve a relaxation of the weighted-cost optimal transport problem with a soft-constraint for the target distribution. In addition, we show that this procedure also applies to the construction of triangular velocity fields via dynamic optimal transport yielding optimal velocity fields. This justifies various variational methodologies for estimating KR maps in practice by minimizing a divergence between the target and pushforward measure through an approximate map. Moreover, it opens the possibilities for novel static and dynamic OT estimators for KR maps.
\end{abstract}
\maketitle

\section{Introduction}
\label{sec:intro}

Transportation of measures seeks a way of transforming one probability measure into another. Among all such transformations between a source measure $\mu$ and target measure $\nu$ defined on $\R^d$, optimal transport seeks the map that also minimizes the average displacement of moving mass from $x$ to $y$ as measured by a cost function $c(x,y)$. For instance, $c$ can be the Euclidean distance squared $c(x,y) = \|x-y\|_2^2$. Alternatively, \cite{carlier2008knothes} considers a class of weighted Euclidean cost function depending on a parameter $\epsilon > 0$ where the weights decay with the index of the variables as follows:
\begin{equation} \label{eq:weighted_cost}
c_{\epsilon}(x, y) = \|x - y\|_\epsilon^2 \coloneqq  \sum_{i=1}^d \eps^{i-1}\left(x_i-y_i\right)^2.
\end{equation}
For a given cost function $c$, source measure $\mu$ and target measure $\nu$, optimal transport then seeks the map $T \colon \R^d \rightarrow \R^d$ that minimizes this cost and satisfies the pushforward constraint $T_{\sharp}\mu = \nu$, which is defined by $\nu(A) = \mu(T^{-1}(A))$ for all Borel sets $A \subset \R^d$. 
The \emph{Monge formulation of optimal transport} for the class of cost functions $c_\eps(x,y)$ is given by
\begin{equation} \label{eq:hard_constraint_Monge} 
    T_\epsilon = \argmin_{T} \left\{  \int_{\mathbb{R}^d} c_\eps(x,T(x)) d\mu : T_\sharp \mu 
    = \nu \right\}.
\end{equation}
Here, $T_\epsilon$ denotes the optimal map for transport cost $c_\epsilon(x,y)$ if it exists, which is nicely characterized by Brenier's theorem for $\epsilon>0$ guaranteeing that the optimal map is the gradient of a convex function \cite{Brenier1991PolarFA}. 

In the limit $\epsilon \rightarrow 0$, \cite{carlier2008knothes} %
showed that the solution of~\eqref{eq:hard_constraint_Monge} converges %
to a map with a \emph{triangular} structure, known as the Knothe-Rosenblatt (KR) transport introduced in~\cite{knoth, rosenblatt}. For measures $\mu,\nu$ defined on $\R^d$, the KR map $T_{\KR}\colon \R^d\to\R^d$ pushing forward $\mu$ to $\nu$ is defined as the unique lower-triangular transformation
\begin{equation}\label{eq:TKR}
    T_{\KR}(x)=(T_1(x_1),T_2(x_1,x_2),\ldots,T_d(x_1,x_2,\dots, x_d)),
\end{equation}
where the $i$-th component $x_i \mapsto T_i(x_1,\ldots,x_{i-1},x_i)$ corresponds to the monotone rearrangement pushing forward the one-dimensional marginal conditional $\mu_i(dx_i|x_1,\dots,x_{i-1})$ to $\nu_i(dx_i|x_1,\dots,x_{i-1})$ for each $(x_1,\dots,x_{i-1}) \in \R^{i-1}$ and $i=1,\dots,d$. Note that the lower-triangular structure of the KR map implies its Jacobian is not symmetric. This contrasts with Brenier maps whose Jacobians are symmetric as they are Hessians of convex scalar functions. Thanks to the structure of its components, the KR map has been used to construct conditional sampling methods with applications in uncertainty quantification \cite{Marzouk_2016}, Bayesian statistics \cite{hosseini2024conditionaloptimaltransportfunction, wang2023efficient}, and generative modeling~\cite{papamakarios2021normalizing, baptista2023conditionalsamplingmonotonegans}. Computationally, KR maps are often identified by searching over a parametrized subset of functions that have the monotone and lower triangular form in~\eqref{eq:TKR}; see~\cite{Baptista_2023, zech2022sparse, irons2022triangular, wang2022minimax} for the optimization, approximation and statistical estimation of KR maps.

However, enforcing the push-forward constraint exactly is difficult in practice as it is generally nonlinear in the transport map. One way to approximate the constraint is to relax it by introducing a probability divergence $\mathcal{D}: \mathcal{P}(\R^d) \times \mathcal{P}(\R^d) \rightarrow \R_+$ that evaluates the discrepancy between the pushforward measure $T_\sharp \mu$ and the target measure $\nu$. A frequently used choice is the Kullback-Leibler (KL) divergence given by
$\dkl(\nu_1||\nu_2) = \int \log(\frac{d\nu_1}{d\nu_2}) d\nu_1$. This yields a family of \emph{soft-constrained Monge problems} indexed by a regularization parameter $\lambda \geq 0$ and cost parameter $\epsilon > 0$ given by
\begin{equation} \label{eq:soft_constraint_Monge} 
    T_{\epsilon,\lambda} 
= \argmin_{T}\left\{\lambda \dkl\left(T_\sharp \mu \|\nu \right) + \int_{\mathbb{R}^d} c_\eps(x,T(x)) d\mu\right\}.
\end{equation}
We denote the family of optimal transport maps that solve each soft-constrained Monge problem as $(T_{\epsilon,\lambda})$. Naturally, one may ask if in the limit $\epsilon \rightarrow 0$ for each $\lambda,$ the solution $T_{\epsilon,\lambda}$ to \eqref{eq:soft_constraint_Monge} also converges to a Knothe-Rosenblatt map as in the hard-constraint setting \eqref{eq:hard_constraint_Monge}. In \cref{thm3}, we prove that this is in fact true. Specifically, we show that for each $\lambda >0$, the solution $T_{\epsilon,\lambda}$ converges to a KR map $T_{\KR,\lambda}$ pushing forward $\mu$ to an intermediate measure $\nu_\lambda$ depending on the regularization parameter $\lambda$. Moreover, we prove in \cref{thm4} that in the limit $\lambda \rightarrow \infty,$ the measure $\nu_\lambda$ and the map $T_{\KR,\lambda}$ converges to the original target measure $\nu$ and the KR map $T_{\KR}$ between $\mu$ and $\nu$, respectively. Finally, in \cref{thm2} we show that this convergence also occurs by taking the limits in the other order, i.e.,  in the limit $\lambda \to \infty$ and then $\epsilon \to 0$, $\{T_{\epsilon,\lambda} \}$ converges to $T_{\KR}$ the Knothe-Rosenblatt map that satisfies the pushforward constraint exactly, i.e., $(T_{\KR})_\sharp\mu = \nu$. This switching of limits is of interest as it allows us to guarantee the existence of a diagonal sequence when sending $(\eps,\lambda)\to (0,\infty)$, where we have that $T_{\eps,\lambda}\to T_\KR$. This is important because practitioners often use algorithms where the variables $(\eps,\lambda)$ are coupled and taken to their limits 0 and $\infty$ simultaneously; see~\cite{baptista2024conditionalsimulationentropicoptimal} for one example. Our result justifies this choice. We also emphasize that our result shows that the optimal transport map solving the soft constrained problem (i) is given by the gradient of a convex function and (ii) converges to a triangular map which cannot be the gradient of a scalar function. In other words, we complete the commutative diagram in \cref{fig:commdiagram}.
\begin{figure}[!ht]
\begin{tikzpicture}[node distance=3cm, every node/.style={anchor=center}]
    \node (T_e_lambda) {$T_{\epsilon,\lambda}$};
    \node (T_e) [right of=T_e_lambda] {$T_\epsilon$};
    \node (T_K_n) [below of=T_e_lambda] {$T_{\KR,\lambda}$};
    \node (T_K) [right of=T_K_n] {$T_{\KR}$};

    \draw[->] (T_e_lambda) -- (T_e) node[midway, above] {$\lambda\to\infty$};
    \draw[->] (T_e_lambda) -- (T_K_n) node[midway, left] {$\epsilon\to 0$};
    \draw[->] (T_e) -- (T_K) node[midway, right] {$\epsilon\to 0 \text{ (see \cite{carlier2008knothes})}$};
    \draw[->] (T_K_n) -- (T_K) node[midway, below] {$\lambda\to\infty$};
\end{tikzpicture}
\caption{``Commutative diagram'' of limits between the minimizers of various optimal transport problems we consider in this work.}
\label{fig:commdiagram}
\end{figure}

The convergence results are stated formally in terms of the narrow convergence of the transport plans that solve Kantorovich problems, which are relaxations of~\eqref{eq:hard_constraint_Monge} and~\eqref{eq:soft_constraint_Monge} and will be introduced in \cref{sec:background}. The convergence of the corresponding transport maps, when they exist, arises as a corollary. As part of our analysis, we show a  general stability result in \cref{thm:kr stab} for the limiting plan that solves the Kantorovich problem with the cost $c_\epsilon$ and varying source and target measures $\mu_\epsilon \rightharpoonup \mu$ and $\nu_\epsilon \rightharpoonup \nu$. In particular, we show that in the limit $\epsilon \rightarrow 0$, the limit of the optimal plan between $\mu_\eps$ and $\nu_\eps$ is also induced by a KR map between $\mu$ and $\nu$ as long as the source and target measures converge to $\mu$ and $\nu$ respectively. \cref{thm4} then follows as a special case when the source measure is fixed and the target measure is changing with the regularization parameter $\lambda$.

Another way of viewing optimal transport is through the Benamou-Brenier \cite{Benamou2000ACF} formulation, also referred to as dynamic optimal transport. Let the source and target measures have density with respect to the Lebesgue measure on $\R^d$, i.e., $\mu(dx)=\rho_0(x)dx$ and $\nu(dx)=\rho_1(x)dx$. Let $v\colon \R^d \times [0,1] \to\R^d$ be a time-varying velocity field, and let $X(t,x)$ be a flow starting from the initial condition $x \in \R^d$ drawn from $\mu$ that satisfies
\begin{equation} \label{eq:flowmap}
    \begin{cases}
       \frac{\partial}{\partial t}X(t,x)=v(t,X(t,x))\\
        X(0,x)=x.
    \end{cases}
\end{equation}
We denote the marginal density of the flow as the push-forward of the source density, $\rho(\cdot,t) = X(t,\cdot)_{\sharp}\rho_0$. Observe that $\rho$ weakly solves the continuity equation
\begin{equation}\label{CE}
    \frac{\partial}{\partial t}\rho(t,x) + \mathrm{div}(v(t,x)\rho(t,x)) = 0.
\end{equation}
Given a pair $(\rho,v)$ that satisfies the continuity equation, we have the Benamou-Brenier formulation for cost function $c_\epsilon(x,y) = \|x-y \|_\epsilon^2$
\begin{align}\label{eq:bb form}
    (\rho_\epsilon,v_\epsilon) = \argmin_{\rho,v} &\left\{\int_0^1\int_{\R^d}\|v(x,t)\|_\epsilon^2 \,\rho_t(x)\ dxdt \right.\\
    &\quad \left. \text{s.t. } \frac{\partial \rho}{\partial t} + \mathrm{div}(v\rho) = 0,\, \rho(0,x) = \rho_0(x), \rho(1,x) = \rho_ 1(x) \right\}. \nonumber
\end{align}
In \cref{thm1} we show that the minimizing velocity field $v_\epsilon$ converges to a field $v_{\KR}$ induced by the Knothe-Rosenblatt map, i.e., the flow map $X(t,x)$ solving ~\eqref{eq:flowmap} with $v = v_{\KR}$ defines a triangular transport pushing forward $\mu$ to $\nu$ by setting $T_{\KR}(x)=X(1,x)$. Furthermore, we consider a soft-constrained version of \eqref{eq:bb form}
\begin{align}\label{eq:softbb}
    (\rho_\epsilon,v_\epsilon) = \argmin_{\rho,v} &\left\{ \lambda \mathcal{D}_{KL}(\rho(1,\cdot)\| \rho_1(\cdot)) + \int_0^1\int_{\R^d}\|v(x,t)\|_\epsilon^2 \,\rho_t(x)\ dxdt \right.\\
    &\quad \left. \text{s.t. } \frac{\partial \rho}{\partial t} + \mathrm{div}(v\rho) = 0,\, \rho(0,x) = \rho_0(x)\right\}. \nonumber
\end{align}
For the soft-constraint dynamic formulation \eqref{eq:softbb}, it is not clear whether the same commutation of limits $\epsilon\to 0$ and $\lambda\to\infty$ holds as it does for the optimal transport plans. We leave this question as an open problem. 

The theoretical questions we examine in this work are related to the topic of unbalanced optimal transport (UOT)~\cite{liero2018optimal,chizat2018interpolating, sejourne2023unbalanced}. UOT consider relaxations of the Kantorovich Problem (see \cref{sec:background}) in which satisfying the source and target constraints exactly is replaced by minimizing divergences between probability distributions. Our work can be seen as a special case of UOT when only the target constraint is relaxed and the cost is given by the weighted Euclidean norm in~\eqref{eq:weighted_cost}. At the same time, this work \emph{extends} the results shown in~\cite{carlier2008knothes} to UOT for the convergence of the minimizer to~\eqref{eq:soft_constraint_Monge}. As far as we know, this is the first result that shows there is a continuation of the solution for the UOT problem between non-triangular (i.e., Brenier maps with $\epsilon =  1$) and triangular maps (i.e., $\epsilon \rightarrow 0$). This result is insightful because it shows that the set of Brenier maps (i.e., gradients of convex functions with symmetric Jacobians) is not closed as we found a family of sequences for different UOT objectives that all converge to triangular maps with non-symmetric Jacobians.

Optimization problems of the form in~\eqref{eq:soft_constraint_Monge} and~\eqref{eq:softbb} have been successfully applied to approximate transports for sampling and conditional generative modeling. For example, \cite{wang2023efficient} seeks a transport map within a parametrized family of (block)-triangular transport maps whose components push forward a known reference to marginal conditional distributions of the target. These components are chosen individually to minimize the KL divergence to match the pushforward constraint while being optimal with respect to the optimal transport problem with the standard Euclidean norm ($\epsilon = 1$). The results of this paper show that the same set of minimizers for the components is attained as the limit of a sequence of UOT problems \eqref{eq:soft_constraint_Monge}
for the entire map as $\epsilon \to 0$ and $\lambda\to\infty$. Similar results are expected to hold for the corresponding dynamic formulation.

While the main framing of this paper regularizes the OT problem \eqref{eq:hard_constraint_Monge} by relaxing the terminal constraint with a KL divergence, this can also be interpreted as regularization of the KL minimization problem through optimal transport. When the source and target measures do not have the same support, the KL divergence may be infinite, but the Wasserstein regularization transports mass from the source to an intermediate target measure such that the KL divergence between the intermediate and target measures is finite. This combination of KL divergence of optimal transport costs is known as the Wasserstein proximal regularization of KL divergence \cite{ambrosio2008gradient,salim2020wasserstein,parikh2014proximal,gu2024combining,baptista2025proximal}. Therefore, \eqref{eq:soft_constraint_Monge} can be interpreted as a Wasserstein proximal regularization for a weighted transport cost. Moreover, this provides a justification for seeking monotone triangular transport maps that minimize KL divergence wherein this particular space of function is induced by a sequence of Wasserstein spaces with the limit of the weighted Euclidean costs.

We hope the results of this paper inspire future computational tools for approximating KR maps as the limit of a sequence of Brenier maps that solve relaxed optimal transport problems. A property of KR maps is that it requires a choice of variable ordering. In practice, each ordering results in a different KR rearrangement and a different approximation error, and thus must be chosen judiciously in advance. This perspective suggests a learning strategy for choosing a variable ordering by adapting the scalings in the weighted Euclidean norm for different variables. Moreover, the particular scaling of the weights $\epsilon$ can be chosen to learn or optimize the orderings in the KR map that yields the most parsimonious representations of the map.

\subsection{Organization} We begin in \cref{sec:background} by reviewing relevant notions in optimal transport, including the Kantorovich formulation, and triangular transport maps. In Section \ref{section2}, we show the first set of results pertaining to Knothe-Rosenblatt maps arising from the soft-constrained OT problem. In \cref{section3}, we show analogous results for the dynamic formulation of the soft constrained problem. In \cref{sec:discussion} we discuss the implications of our results for computation and potential directions for future work. Detailed proofs and supporting auxiliary results are postponed to \cref{app:main_proofs,app:additional_proofs}.

\section{Background, Notation and Assumptions}
\label{sec:background}

\subsection{Formulations of Optimal Transport} We start by recalling three canonical problems of optimal transportation of measures. Given a source measure $\mu$ and target measure $\nu$, optimal transport considers the best way to move mass from one probability distribution to another. In static optimal transport, the transport is found by minimizing the average cost of displacing mass with respect to a cost function $c\colon \R^d\times\R^d\to \R$. The Monge problem is formulated as
\begin{align}
    M(\mu,\nu) \coloneqq \inf_T\left\{\int_{\R^d}c(x,T(x))\ d\mu(x):T_{\sharp}\mu=\nu\right\},
\end{align}
where we want to find the transport map $T\colon\R^d\to \R^d$ that pushes forward mass from probability distribution $\mu$ to $\nu$ such that the transport cost is minimized. Importantly, $T$ is Borel measurable. Often times, finding a transport maps that satisfies the pushforward constraints is difficult, or may not be possible (e.g., consider $\mu = \delta_{x_0}$, $\nu = \frac{1}{2}\delta_{y_0}+ \frac{1}{2}\delta_{y_1}$, where $y_0 \neq y_1$).

Instead, a relaxation of the Monge problem seeks a joint distribution $\gamma \in \mathcal{P}(\R^d \times \R^d)$, also referred to as a coupling,  %
with marginal distributions $\mu$ and $\nu$, i.e., %
$(\pi_X)_{\sharp}\gamma=\mu$ and $(\pi_Y)_{\sharp}\gamma=\nu$ where $\pi_X(x,y)=x$ and $\pi_Y(x,y)=y$ are projection operators. We denote the space of all such couplings as $\Pi(\mu,\nu)$; the Kantorovich problem is given by 
\begin{align} \label{eq:kantorovichproblem}
    K(\mu,\nu):=\inf_\gamma \left\{\int_{\R^d\times\R^d}c(x,y)\ d\gamma(x,y):\gamma\in\Pi(\mu,\nu)\right\}.
\end{align}
Solutions of the \eqref{eq:kantorovichproblem} are referred to as \emph{optimal plans}. Under mild conditions, namely the cost function $c$ being lower semi-continuous and nonnegative, the Kantorovich problem is guaranteed to have a minimizer; see~\cite{otam, figalli2021invitation} for more details. When a transport map exists, the value of the Monge and Kantorovich problems coincide. The following theorem presents conditions where this occurs for the (weighted) quadratic cost as well as characterizing the structure of the optimal map and associated plan. Note that these conditions are not the most general ones one could state, but we will use them here for simplicity of exposition. 
\begin{theorem}[Brenier's theorem \cite{Brenier1991PolarFA}]\label{thm:brenier's}
    Let $c(x, y)=\|x-y\|^2$ (or equivalently $c_\eps(x, y)=\|x-y\|_\eps^2)$. Suppose that $\mu$ and $\nu$ have finite second moment, and that $\mu \ll d x$. Then there exists a unique optimal plan $\gamma$. In addition, $\gamma=(id \times T)_{\sharp} \mu$ and $T=\nabla \varphi$ for some convex function $\varphi$. Further, $T$ is the unique optimal transport map pushing $\mu$ onto $\nu$.
\end{theorem}

\subsection{Triangular Maps and Assumptions.} 
In this subsection, we will outline two construction of the triangular transport maps known as the Knothe-Rosenblatt rearrangement that was introduced in \cref{sec:intro}. The first provides explicit expressions for the components of the KR map based on the marginal conditional cumulative distribution functions.  %
Let us consider measures defined on $\R^2$, the construction for higher dimensions follows analogously. In this case, the KR map acts as follows
\begin{equation*}
\begin{dcases}
    T_{\KR}(x_1,x_2) = (T_1(x_1),T_2(x_1,x_2)),\\
    T_{1\sharp}(f_1 dx_1) = g_1 dy_1,\\
    T_2(x_1,\cdot)_{\sharp}\left(\frac{f(x_1,\cdot)}{f_1(x_1)}dx_2\right)=\frac{g(T_1(x_1),\cdot)}{g_1(T_1(x_1))}dy_2,
\end{dcases}
\end{equation*}
where 
\begin{equation*}
    \begin{dcases}
        d\mu(x_1,x_2) = f(x_1,x_2)dx_1dx_2,\qquad f_1(x_1) = \int_\R f(x_1,x_2) dx_2,\\
        d\nu(y_1,y_2) = g(y_1,y_2)dy_1dy_2,\qquad g_1(y_1) = \int_\R g(y_1,y_2)dy_2.
    \end{dcases}
\end{equation*} 
The components of such a KR map satisfying the above conditions is uniquely given by
\begin{align}
    &T_1(x_1) = (G_1^{-1} \circ F_1)(x_1), \\ 
    &T_2(x_1,x_2) = \left(G_{2|1} ^{-1}(\cdot, x_1) \circ F_{2|1} (\cdot, x_1) \right)(x_2),
\end{align}
where the functions $F_1, G_1$ $F_{2|1}, G_{2|1}$ are cumulative distribution functions
\begin{align}
    &F_1(x_1) = \int_{-\infty}^{x_1} f_1(t) dt, \quad G_1(y_1) = \int_{-\infty}^{y_1} g_1(t) dt, \\ & F_{2|1}(x_2, x_1) = \int_{-\infty}^{x_2} \frac{f(x_1,t)}{f_1(x_1)} dt\quad G_{2|1}(y_2, y_1) = \int_{-\infty}^{y_2} \frac{g(y_1,t)}{g_1(y_1)} dt.  \nonumber
\end{align}

While this characterization is explicit, it is difficult to compute in practice, especially in high dimensions, because it requires repeated integration of the joint density to extract the marginal conditional densities. The following result from \cite{carlier2008knothes} provides an alternative characterization of the KR map via a sequence of optimal transport problems. We consider a special class of cost functions, but we refer the reader to the original result for a generalization. 

\begin{theorem}[Theorem 2.1 in \cite{carlier2008knothes}]\label{thm0}
Suppose that our source and target measure $\mu$ and $\nu$ are absolutely continuous and have finite second moment. Let $\gamma_{\eps}$ be the optimal transport plan between $\mu$ and $\nu$ for the costs $c_{\varepsilon}(x, y)=\sum_{i=1}^d \eps^{i-1}\left(x_i-y_i\right)^2$. Let $T_{\KR}$ be the Knothe-Rosenblatt map between $\mu$ and $\nu$ and $\gamma_{\KR} \in$ $\mathcal{P}\left(\mathbb{R}^d \times \mathbb{R}^d\right)$ the associated transport plan (i.e. $\left.\gamma_{\KR}:=(i d \times T_{\KR})_\sharp \mu\right)$. Then $\gamma_{\eps} \rightharpoonup \gamma_{\KR}$ as $\eps \rightarrow 0$.

Moreover, the plans $\gamma_{\eps}$ are induced by transport maps $T_{\eps}$ that converge to $T_{\KR}$ in $L^2(\mu)$ as $\eps \rightarrow 0$. 
\end{theorem}

In \cite{carlier2008knothes}, the authors deal with the case of hard-constrained optimal transport. Computing either the KR map or the optimal transport map remains challenging if the pushforward constraint is enforced exactly. In this work we extend this result to unbalanced OT, i.e., the soft-constrained problem, which is solved in practice. To do so, we make the following assumptions on our source and target measures $\mu$ and $\nu$.
\begin{assumption}\label{ass:source_target_measures}
\hspace{100pt}
\begin{itemize}
    \item[{$(\mathcal{A}_1)$}] The source measure $\mu$ and target measure $\nu$ are absolutely continuous with respect to the Lebesgue measure;
    \item[{$(\mathcal{A}_2)$}] The source and target measures have compact support.
\end{itemize}
\end{assumption}
For some of our results, we can relax these assumptions to
\begin{assumption}\label{asm:ac gaussian}
\hspace{100pt}
\begin{itemize}
    \item[{$(\mathcal{B}_1)$}] The source and target measure are absolutely continuous with respect to any fully supported, independent across coordinates (like the Gaussian) probability measure on $\R^d$;
    \item[{$(\mathcal{B}_2)$}] The source and target measure have finite second moment.
\end{itemize}
\end{assumption}
Indeed, most of our results can be generalized to hold under \cref{asm:ac gaussian} instead of \cref{ass:source_target_measures}, see \cref{rmk:assAvsB}. We consider the stronger \cref{ass:source_target_measures} here for simplicity of exposition.

\subsection{Notation}\label{sec:notat} Throughout the paper, we refer to $\mathcal{P}(\R^d)$ as the space of probability measures over $\R^d$ and $\mathcal{P}(\R^d\times\R^d)$ as the space of probability measures of $\R^d\times\R^d$. $\mathcal{P}_p(\R^d)$ refers to probability measures with $p$-th finite moment. If not specified, all the measures we work with are probability measures in $\mathcal{P}(\R^d)$ defined over the Borel sets. We refer to $\Pi(\mu,\nu)$ as the set of couplings (transport plans) in $\mathcal{P}(\R^d\times\R^d)$ with marginals $\mu$ and $\nu$.  $C^k(\R^d)$ refers to the space of functions with $k$ continuous derivatives, and $C_b(\R^d)$ refers to functions that are continuous and bounded. The notation $\rightharpoonup$ refers to narrow convergence of measures (i.e., convergence in duality with $C_b$.) We also talk about the notion of two measures being equal, $\mu =\nu$, when for all $\varphi\in C_b$, $\int \varphi d\mu =\int \varphi d\nu.$ We define the projection operators $\pi_X(x,y)=x,
    \pi_Y(x,y)=y,
    \pi_{i:j}(x_1,\ldots,x_d)=(x_i,\ldots,x_j) \text{ for } 1\leq i\leq j\leq d.$
 Given a coupling $\gamma$ on $\R^d\times \R^d$, define the projection operator $\pi^{x,y}_{k:r}: \R^d\times \R^d \to \R^{r-k+1} \times \R^{r-k+1}$ for $1\le r\le k \le d$ where $\pi^{x,y}_{k:r}(x_1,\ldots x_d, y_1, \ldots, y_d) = (x_k, \ldots, x_r, y_k, \ldots, y_r) \coloneqq (x_{k:r},y_{k:r})$. Analogously, define $\pi_k^{x,y}: \R^d\times \R^d \to \R \times \R$, where $\pi^{x,y}_{k}(x,y) = (x_k,y_k)$, and $\pi_k^x: \R^d \to \R$. Furthermore, given a coupling $\gamma(x,y)$, denote $\gamma^{r:k} = (\pi_{r:k}^{x,y})_\sharp\gamma $, the projection of $\gamma$ only on the $r$ through $k$ marginals of $x$ and $y$. When $r = k$, we omit the colon. For any probability measure $\mu \in \mathcal{P}(\R^d)$, its conditional distribution  $\mu_{x_{r:k}}$ is denoted with subscript notation and given by
 \begin{align*}
     \mu_{x_{r:k}}(dy_{1:r-1}, dy_{k+1:d}) = \frac{\mu(dy_{1:r-1}, x_{r:k}, dy_{k+1:d})}{\int \mu(z_{1:r-1}, x_{r:k},z_{k+1:d})dz_{1:r-1}, dz_{k+1:d}};
 \end{align*}
its marginal distributions $\mu^{r:k}$ is denoted with a superscript notation and given by 
 \begin{align*}
     \mu^{r:k}(dy_{r:k}) = \int \mu(z_{1:r-1}, dy_{r:k},z_{k+1:d})dz_{1:r-1}, dz_{k+1:d}\,.
 \end{align*}
 By the disintegration theorem, we have $\mu = \mu^{r:k}\otimes\mu_{x_{r:k}}$. Conditional-marginal distributions are then similarly defined; for example $\mu^{1:k}_{x_{1:r}}$ for $k>r$ denotes the conditional distribution of $x_{1:k}$ given $x_{1:r}$.

\section{Knothe-Rosenblatt Maps via Soft-constraints}\label{section2}

In this section we show that the soft-constrained Kantorovich problem gives us a Knothe-Rosenblatt plan from $\mu$ to an intermediate measure $\nu_\lambda$ by considering the limit $\eps\to0$ in the weighted Euclidean cost. Moreover, we show that that plan converges to the Knothe-Rosenblatt plan from $\mu$ to $\nu$ in the limit $\lambda\to\infty$. 
\cref{sec:main_results_softconstraint} presents an overview of the main results. \cref{sec:convergence_soft_constraint} details the results on the convergence of the solution of the soft-constrained transport problems as $\epsilon \rightarrow 0$ for any $\lambda>0$. Similarly, \cref{sec:lambdatoinfty} shows the convergence of the optimal maps as both $\epsilon \to 0 $ and $\lambda \to \infty$ in any order.  Complete proofs of the main results are found in \cref{app:main_proofs} and supporting lemmas for the proofs are found in \cref{app:additional_proofs}. Also in \cref{app:additional_proofs}, we show existence and uniqueness of the minimizer for the soft-constrained problem, we show that the minimum of the soft-constrained problem is always bounded above by the minimum of the hard-constraint problem, and we prove a key result that we use in our proofs in this section: stability of conditionals and thus stability of KR plans with respect to weak convergence of the source and target measures.

\subsection{Statement of Main Results} \label{sec:main_results_softconstraint}

We first present the soft-constrained Kantorovich problem with the weighted Euclidean cost. 

\begin{definition}[Soft-constrained Kantorovich problem]
The relaxed Kantorovich cost is defined as
    \begin{align}\label{eq:soft_kantorovich}
    K_{\eps,\lambda}(\mu,\nu) = \min _{\gamma}\left\{ \lambda\dkl\left((\pi_Y)_\sharp \gamma \| \nu\right)+ \int_{\mathbb{R}^d\times\R^d} c_\eps(x,y) d\gamma:(\pi_X)_\sharp\gamma=\mu\right\}.
    \end{align}
\end{definition}
Here, we emphasize that we do not require that the solution satisfies the target marginal constraint --- only that the second marginal is a perturbation of the target measure.  We can also define the corresponding hard-constrained Kantorovich problem
\begin{definition}[Hard-constrained Kantorovich Problem]
\[
K_{\eps}(\mu,\nu)=\min _{\gamma}\left\{\int_{\mathbb{R}^d} c_\eps(x,y) d\gamma:\gamma\in\Pi(\mu,\nu)\right\}.
\]
\end{definition}

\begin{rem}
    In the definition of $K_{\eps,\lambda}(\mu,\nu)$ the choice of $\mathcal{D}_{\KL}$ can be replaced by any member of the class of $\varphi$-divergences. %
    Borrowing from the definition in \cite{sejourne2023unbalanced}, these are divergences that satisfy the following.
    Let $\varphi\colon(0,\infty)\to[0,\infty]$ be an entropy function that is convex and lower semi-continuous such that $\varphi(1) = 0$. Define $\varphi'_\infty := \lim_{p\to\infty}\varphi(p)/p$. For any $\mu,\nu\in \mathcal{P}(\R^d)$, write the Radon-Nikodym decomposition of $\mu = \frac{d\mu}{d\nu}\nu +\mu^{\perp}$. The $\varphi-$divergence is then defined as 
    \[
    \mathcal{D}_{\varphi}(\mu\|\nu):= \int_{\R^d} \varphi\left(\frac{d\mu}{d\nu}\right)d\nu(x) + \varphi'_\infty\int_{\R^d}d\mu^{\perp}.
    \]
\end{rem}

The first set of results deals entirely with static optimal transport. To our knowledge, it is the first attempt to bridge soft-constrained optimal transport with Knothe-Rosenblatt maps. In practice, the convexity introduced to the soft-constrained problem makes it  more tractable to compute, and the following results prove theoretical guarantees for a tractable approach to estimating Knothe-Rosenblatt maps. Our main contributions are as follows:
\begin{enumerate}
    \item \cref{thm3} shows that as $\eps\to0$, the minimizer of $K_{\eps,\lambda}(\mu,\nu)$ is given by a plan $\gamma_{\KR,\lambda}$ which induces a Knothe-Rosenblatt map $T_{\KR,\lambda}$. This triangular map pushes forward $\mu$ to some intermediate measure $\nu_\lambda$ which is a member of a sequence indexed by $\lambda$ that converges to $\nu$ as $\lambda\to\infty$. %
    \item \cref{thm4} shows that taking $\eps\to 0$, and then $\lambda\to\infty$, the minimizer to $K_{\eps,\lambda}(\mu,\nu)$ converges narrowly to the Knothe-Rosenblatt plan from $\mu$ to $\nu$.
    \item \cref{thm2} shows that reversing the order of the limits, i.e., taking $\lambda\to\infty$ and then $\eps\to 0$, the minimizer of $K_{\eps,\lambda}(\mu,\nu)$ also converges narrowly to the Knothe-Rosenblatt plan from $\mu$ to $\nu$.
\end{enumerate}
We emphasize that convergence of the transport plans to a KR plan yields convergence of transport maps to a KR map immediately, as long as the source measure is absolutely continuous with respect to the Lebesgue measure.

\subsection{Convergence to a soft-constrained Knothe-Rosenblatt map} \label{sec:convergence_soft_constraint}

To prove the above results, we also show that the set of KR plans with certain conditions is closed under narrow convergence and that if $\gamma_{\KR,n}$ is a KR plan from $\mu_n$ to $\nu_n$ \footnote{Note that the source measure $\mu_n$ is also changing with $n$, in contrast to \eqref{eq:soft_constraint_Monge}, \eqref{eq:soft_kantorovich}. } with $\mu_n\rightharpoonup\mu$ and $\nu_n\rightharpoonup\nu$, then $\gamma_{\KR,n}\rightharpoonup\gamma_{\KR}$, where $\gamma_{\KR}$ is the KR plan from $\mu$ to $\nu$. To the best of our knowledge, this is the first proof showing KR plans are closed under narrow convergence. Relatedly, in \cite{beiglböck2023knotherosenblattdistanceinducedtopology}, the authors showed that triangular transport maps are closed under $L^p$ convergence.

 \begin{theorem}[From OT to KR]\label{thm:kr stab}
        
        Let our cost be the scaled Euclidean cost $c_\eps(x,y)$ and let $\mu_\eps$ and $\nu_\eps$ be two sequence of probability measures that are uniformly compactly supported with $\mu_\eps\rightharpoonup \mu$ and $\nu_\eps\rightharpoonup\nu$. Further suppose their densities exists and converge pointwise a.e. to the densities of $\mu$ and $\nu$. For each $\eps>0$, let $\gamma_\eps \in \Pi(\mu_\eps,\nu_\eps)$ be an optimal transport plan with respect to $c_\eps(x,y)$. Then $\gamma_\eps \rightharpoonup \gamma_{\KR}$ as $\eps\to 0$ where $\gamma_{\KR}$ is the Knothe-Rosenblatt plan from $\mu$ to $\nu$. If $\mu\ll dx$, it follows that $\gamma_{\KR}$ is induced by a Knothe-Rosenblatt map $T_{\KR}$; and if $\mu_\eps = \mu$ for all $\eps$, then $T_\eps$, the optimal map under $c_\eps(x,y)$ from $\mu$ to $\nu_\eps$, converges to $T_{\KR}$ in $L^2(\mu)$.
    \end{theorem} 

    \begin{rem}
        In this theorem, we only assume that $(\mu_\eps)_\eps$ and $(\nu_\eps)_\eps$ are uniformly compactly supported respectively. Since $\mu_\eps \rightharpoonup \mu$ and $\nu_\eps \rightharpoonup \nu$, it follows that $\mu$ and $\nu$ are also supported in the same compact sets, respectively. Hence, we do not state a compact support assumption for $\mu$ and $\nu$ separately. %
    \end{rem}

    \begin{proof}[Proof Sketch]
        To prove this, we start by observing that 
        \begin{equation}\label{eq:limit 1}
        \int_{\R^d\times\R^d} c_\eps(x,y)d\gamma_\eps \leq \int_{\R^d\times\R^d}c_\eps(x,y)d\gamma_{\KR,\eps}\,,
        \end{equation}
        where $\gamma_{\KR,\eps}$ is the Knothe-Rosenblatt plan between $\mu_\eps$ and $\nu_\eps$. Note that $\gamma_\eps\rightharpoonup\gamma$ by stability of optimal transport plans (\cref{thm stab}), and stability of KR maps (\cref{lem:closed KR}) tells us that $\gamma_{\KR,\eps}\rightharpoonup\gamma_{\KR}$. By assumption, $\gamma_\eps$ and $\gamma_{\KR,\eps}$ are uniformly compactly supported; on this support, $c_\eps$ acts as a test function in $C_b(\R^d\times\R^d)$, and so passing to the limit $\eps\to 0$ in \eqref{eq:limit 1} gives
        \begin{equation}\label{eq:limit 2}
        \int_{\R^d\times\R^d} (x_1-y_1)^2 d\gamma \leq \int_{\R^d\times\R^d} (x_1-y_1)^2 d\gamma_{\KR}\,.
        \end{equation}
        Since $\gamma_{\KR}$ is quadratically optimal on its first coordinate, this tells us that $\gamma_{\KR}$ and $\gamma$ agree on its first marginals. Repeating a similar strategy to \cite{carlier2008knothes}, we can show by induction that $\gamma_{\KR}$ and $\gamma$ agree on all its marginals. Thus $\gamma=\gamma_{\KR}$. The full details for this proof are deferred to \cref{app:main_proofs}.
    \end{proof}

This is, in spirit, the exact same result as \cite{carlier2008knothes}, \cite{otam} with crucial difference that here the source and target measures are dependent on the cost scaling parameter $\epsilon$. To show this, we prove that the KR plans with certain conditions are closed under narrow convergence and that all marginal conditionals of these KR plans also narrowly converge; see proof of \cref{lem:closed KR}.

\begin{theorem}[Convergence to a perturbed Knothe-Rosenblatt map]\label{thm3}
Assume the source $\mu$ and target $\nu$ satisfy \cref{ass:source_target_measures}. Let $\gamma_{\eps,\lambda}$ be the minimizer to $K_{\eps,\lambda}(\mu,\nu)$. As we take $\eps\to 0$, we get $\gamma_{\eps,\lambda} \rightharpoonup \gamma_{\KR,\lambda}$ where $\gamma_{\KR,\lambda}$ is Knothe-Rosenblatt plan (with associated Knothe-Rosenblatt map $T_{\KR,\lambda}$) from $\mu$ to some target $\nu_\lambda$.
\end{theorem}

\begin{proof}
    To prove this theorem, we need two ingredients. Firstly, observe by \cref{cor:unique min} that $\gamma_{\eps,\lambda}$ is also the optimal coupling in the hard constraint problem from $\mu$ to $\nu_{\eps,\lambda}:=(\pi_Y)_\sharp\gamma_{\eps,\lambda}$. That is, \[
    \gamma_{\eps,\lambda} = \argmin_\gamma \left\{\int_{\R^d\times\R^d}c_\eps(x,y):\gamma\in\Pi(\mu,\nu_{\eps,\lambda})\right\}.
    \] Thus, it suffices to apply \cref{thm:kr stab} to $\gamma_{\eps,\lambda}$, for which we need the following to hold:
\begin{enumerate}
    \item $\nu_{\eps,\lambda}:= (\pi_Y)_\sharp\gamma_{\eps,\lambda}$ has finite second moment and has compact support;
    \item $\nu_{\eps,\lambda}$, as $\eps\to 0$, has a weak limit  we call $\nu_\lambda$;
    \item The density of $\nu_{\eps,\lambda}$ exists and converges pointwise to the density of $\nu_\lambda$, which also exists.
\end{enumerate}
In what follows, we carefully verify that each of these conditions is satisfied.
\begin{enumerate}
    \item If $K_{\eps,\lambda}(\mu,\nu)<\infty$, then any minimizer $\gamma_{\eps,\lambda}$ of \eqref{eq:soft_kantorovich} has to satisfy $\nu_{\eps,\lambda}\ll \nu$ which has compact support. This tells us that $\nu_{\eps,\lambda}$ are uniformly compactly supported, from which we deduce that $(\nu_{\eps,\lambda})_{\eps,\lambda}$ has uniformly bounded second moments.
    \item Thanks to (1), the family $(\nu_{\eps,\lambda})_\eps$ is tight. Hence, up to a subsequence, there exists  a probability measure $\nu_\lambda$ such that $\nu_{\eps,\lambda}\rightharpoonup\nu_\lambda$ by Prokhorov's theorem.
    \item Since $\mathcal{D}_{KL}(\nu_{\eps,\lambda}\| \nu)<\infty$ and $\nu\ll dy$, we conclude that $\nu_{\eps,
\lambda}$ is absolutely continuous, and we denote $d\nu_{\eps,\lambda}(y) =g_{\eps,\lambda}(y)dy$ and $d\nu(y)= f(y)dy$. To characterize $\nu_{\eps,\lambda}$, we write the Euler-Lagrange equation of the functional
       \begin{align*}
     \mathcal{L}[\gamma,\varphi] := \lambda\dkl((\pi_Y)_\sharp\gamma\|\nu)+\int_{\R^d\times\R^d} c_\eps(x,y)d\gamma +\int_{\R^d\times \R^d}\varphi(x) d(\gamma-\mu),
     \end{align*} for $\gamma\in\mathcal{P}(\R^d\times\R^d)$ and $\varphi \in C_b(\R^d)$. Note that
     \begin{align*}
    &\gamma_{\eps,\lambda}= \argmin_\gamma\sup_{\varphi\in C_b(\R^d)}\mathcal{L}[\gamma,\varphi]\,,\qquad 
    \varphi_{\eps,\lambda} = \argmax_{\varphi\in C_b(\R^d)} \mathcal{L}[\gamma_{\eps,\lambda},\varphi]\,.
     \end{align*}
     Taking the first variation with respect to $\gamma$ and recalling $\nu_{\eps,\lambda}:= (\pi_Y)_\sharp\gamma_{\eps,\lambda}$, we have 
     \begin{align*}
         \delta_\gamma \mathcal{L}[\gamma,\varphi]\bigg|_{\gamma = \gamma_{\eps,\lambda}\atop\varphi = \varphi_{\eps,\lambda}} &= \delta_\gamma \left(\lambda\dkl((\pi_Y)_\sharp\gamma\|\nu)+\int_{\R^d} c_\eps(x,y)d\gamma + \int_{\R^d}\varphi(x) d((\pi_X)_\sharp\gamma-\mu)\right)\bigg|_{\gamma = \gamma_{\eps,\lambda}\atop\varphi = \varphi_{\eps,\lambda}}\\
         & = \lambda\left(\log\left(\frac{g_{\eps,\lambda}}{f}\right)-\int_{\R^d}\log\left(\frac{g_{\eps,\lambda}}{f}\right)g_{\eps,\lambda}dy\right) + c_{\eps}(x,y)+\varphi_{\eps,\lambda}(x).
     \end{align*}
     The first variation with respect to $\varphi$ evaluated at the Lagrange multiplier corresponding to $\gamma_{\eps,\lambda}$ is 
     \[
     \delta_\varphi \mathcal{L}[\gamma,\varphi]\bigg|_{\gamma = \gamma_{\eps,\lambda}\atop\varphi = \varphi_{\eps,\lambda}} = (\pi_X)_\sharp\gamma_{\eps,\lambda}-\mu.
     \]
     Hence, the Euler-Lagrange conditions are
     \[
     \begin{dcases}
         \lambda\left(\log\left(\frac{g_{\eps,\lambda}(y)}{f(y)}\right)-\int_{\R^d}\log\left(\frac{g_{\eps,\lambda}}{f}\right)g_{\eps,\lambda}dy\right) + c_{\eps}(x,y)+\varphi_{\eps,\lambda}(x) = \mathrm{constant}(\eps,\lambda)\\
          (\pi_X)_\sharp\gamma_{\eps,
          \lambda}(x) = \mu(x)
     \end{dcases}
     \]
     for any $(x,y)\in \mathrm{supp}(\gamma_{\eps,\lambda})$. This implies that there exists $\psi(y)_{\eps,\lambda} := c_\eps(x,y)+\varphi_{\eps,\lambda}(x)$ for $(x,y)\in \mathrm{supp}(\gamma_{\eps,\lambda})$, a complementary slackness condition. Solving for $\nu_{\eps,\lambda}$, we obtain
     \[
     g_{\eps,\lambda}(y)= \frac{1}{Z_{\eps,\lambda}}\exp\left(\frac{-c_\eps(x,y)-\varphi_{\eps,\lambda}(x)}{\lambda}\right)f(y).
     \]
     Here, $Z_{\eps,\lambda}$ is a normalizing constant (independent of $(x,y)\in\mathrm{supp}(\gamma_{\eps,\lambda})$) such that
     \[
     Z_{\eps,\lambda} = \int_{\R^d}\exp\left(\frac{-c_\eps(x,y)-\varphi_{\eps,\lambda}(x)}{\lambda}\right)d\nu(y).
     \]
     Therefore, we can write $\gamma_{\eps,\lambda}-$a.e.
     \begin{align*}
              g_{\eps,\lambda}(y) & = \frac{1}{D_{\eps,\lambda}(x)}\exp\left(\frac{-c_{\eps}(x,y)}{\lambda}\right)f(y)
     \end{align*}  
    where 
     \[
     D_{\eps,\lambda}(x) = \int_{\R^d}\exp\left(\frac{-c_\eps(x,y)}{\lambda}\right)d\nu(y) = Z_{\eps,\lambda}\exp\left(\frac1\lambda\varphi_{\eps,\lambda}(x)\right).
     \]
     Even though $D_{\eps,\lambda}(x)$ depends on $x$, complementary slackness ensures that $\frac{1}{D_{\eps,\lambda}(x)}\exp\left(\frac{-c_\eps(x,y)}{\lambda}\right)$ remains constant in $x$ for fixed $y$ in $\mathrm{supp}(\gamma_{\eps,\lambda}).$ Since $c_\eps(x,y)\to c_0(x,y)=|x_1-y_1|^2$ pointwise, and $\exp(\frac{-c_\eps(x,y)}{\lambda})\leq 1 \in L^1(\nu)$, we can apply dominated convergence to say $D_{\eps,\lambda}(x)\to D_{0,\lambda}(x)$ pointwise,
     where
     \begin{align*}
         D_{0,\lambda}(x)= \int_{\R^d}\exp\left(\frac{-|x_1-y_1|^2}{\lambda}\right)d\nu(y)\,.
     \end{align*}
     Hence, the product $D_{\eps,\lambda}(x)^{-1}\exp(\frac{-c_\eps(x,y)}{\lambda})$ converges pointwise, and by Scheff\'e's Lemma, we have $g_{\eps,\lambda}(y)dy\rightharpoonup g_\lambda(y)dy$, where
     \begin{align*}
       g_\lambda(y):=  \frac{1}{D_{0,\lambda}(x)} \exp\left(\frac{-c_{0}(x,y)}{\lambda}\right)f(y)\,.
     \end{align*}
    By uniqueness of weak limits, we conclude $ \nu_{\lambda}(dy)=g_{\lambda}(y)dy$.
    Now we can apply \cref{thm:kr stab} to conclude that sending $\eps \to 0$ results in $\gamma_{\eps,\lambda}\rightharpoonup \gamma_{\KR,\lambda}$ which is a Knothe-Rosenblatt plan in $\Pi(\mu,\nu_\lambda)$.
\end{enumerate}
    \end{proof}

\subsection{Convergence to the hard-constraint Knothe-Rosenblatt map} \label{sec:lambdatoinfty}
We can also consider $K_{\eps,\lambda}(\mu,\nu)$ but with $\lambda$ not being fixed. What can we say about the minimizer to $K_{\eps,\lambda}(\mu,\nu)$ when we take $\lambda\to \infty$? Does it converge to the Knothe-Rosenblatt map of the hard-constraint problem? We now show that the minimizer of the soft-constraint problem converges to the minimizer of the hard-constraint problem, as $\eps\to 0$ and $\lambda\to\infty$.

\begin{theorem}\label{thm4}
    Let $\gamma_{\eps,\lambda}$ be the minimizer of $K_{\eps,\lambda}(\mu,\nu)$. Then under the same assumptions as in \cref{thm3}, we have that as $\eps\to 0$, $\gamma_{\eps,\lambda}\rightharpoonup \gamma_{\KR,\lambda}$, the Knothe-Rosenblatt plan from $\mu$ to $\nu_\lambda:= (\pi_Y)_\sharp\gamma_{\KR,\lambda}$, and as $\lambda\to\infty$, $\gamma_{\KR,\lambda}\rightharpoonup \gamma_{\KR}$, the Knothe-Rosenblatt plan from $\mu$ to $\nu$.
\end{theorem}
\begin{proof}
Denote $\gamma_\eps:=\arg\min_\gamma\{\int_{\R^d}c_\eps(x,y)d\gamma:\gamma\in\Pi(\mu,\nu)\}$. By \cref{thm3}, as we take $\eps \to 0$, $\gamma_{\eps,\lambda}\rightharpoonup \gamma_{\KR,\lambda}$ narrowly (up to subsequence). Dividing by $\lambda$, we have the following evaluation:
    \begin{align}\label{eq:idk?}
    &\mathcal{D}_{KL}((\pi_Y)_\sharp\gamma_{\KR,\lambda}\|\nu) + \lim_{\eps\to0}\frac{1}{\lambda}\int_{\mathbb{R}^d\times\R^d} c_\eps(x,y) d\gamma_{\eps,\lambda}\notag\\
    &\leq\lim_{\eps\to 0}\left[\mathcal{D}_{KL}((\pi_Y)_\sharp\gamma_{\eps,\lambda}\|\nu) + \frac{1}{\lambda}\int_{\mathbb{R}^d\times\R^d} c_\eps(x,y) d\gamma_{\eps,\lambda}\right]\notag\\
    &\leq\lim_{\eps\to 0}\frac{1}{\lambda}\int_{\R^d}c_\eps(x,y)\ d\gamma_\eps\leq \frac{2M}{\lambda}.
    \end{align}
    The first inequality is from the weak-continuity of KL divergence which we can apply because $(\pi_Y)_\sharp\gamma_{\eps,\lambda}\rightharpoonup(\pi_Y)_\sharp\gamma_{\KR,\lambda}$ when $\eps\to 0$. The second inequality is by \cref{lem:soft-hard}: the soft-constraint minimal value is less than the hard-constraint minimal value. In the last line, we use the fact that there is a uniform bound independent of $\eps$ on $\int_{\R^d}c_\eps(x,y)d\gamma_\eps$. To see this, 
    recall that $\mu, \nu$ have finite second moments, and so there exists a constant $M>0$ such that
    \[
    \int_{\R^d\times\R^d}\left(|x|^2 +|y|^2
    \right) d\gamma_\eps = \int_{\R^d}|x|^2 d\mu+\int_{\R^d}|y|^2d\nu \leq M\,.
    \] Thus, for $\eps\in(0,1)$, 
    \begin{align*}
    \int_{\R^d\times\R^d} c_\eps(x,y)d\gamma_\eps&\leq\int_{\R^d\times\R^d} \|x-y\|^2 d\gamma_\eps
    \leq\int_{\R^d}2|x|^2d\mu+ \int_{\R^d}2|y|^2d\nu \leq 2M.
    \end{align*} Then sending $\lambda\to\infty$ in \eqref{eq:idk?} implies that 
    \[
    \lim_{\lambda\to\infty}\mathcal{D}_{KL}((\pi_Y)_\sharp\gamma_{\KR,\lambda}\|\nu)=0.
    \]
    This implies that $\nu_\lambda=(\pi_Y)_\sharp\gamma_{\KR,\lambda}\rightharpoonup\nu$ narrowly, since $\nu$ is a probability measure. Also, as in the proof of \cref{thm3}, we have that $d\nu_\lambda(y) = g_\lambda(y)dy$ where
    \[
    g_\lambda(y) = \frac{1}{D_{0,\lambda}(x)}\exp\left(\frac{-(x_1-y_1)^2}{\lambda}\right)f(y),
    \]
    and $d\nu(y) = f(y)dy$. To take the limit $\lambda\to\infty$, note that $\exp\left(\frac{-(x_1-y_1)^2}{\lambda}\right)$ converges pointwise to 1, is upper bounded by 1 and is in $L^1(\nu).$ By dominated convergence, $D_{0,\lambda}(x)\to 1$ pointwise as $\lambda\to \infty$. Hence, $\lim_{\lambda\to\infty}g_\lambda(y) = f(y)$.
    Thus, we have enough to conclude that $\gamma_{\KR,\lambda}\rightharpoonup \gamma_{\KR}\in\Pi(\mu,\nu)$ up to subsequence by the closedness of Knothe-Rosenblatt plans, \cref{lem:closed KR}. 
\end{proof}

We now show that changing the order of the limits does not change the minimizer. That is taking $\lambda\to \infty$ first, then $\eps\to 0$ still results in the soft-constraint minimizer converging to the same Knothe-Rosenblatt plan.

\begin{theorem}\label{thm2}
     Assume the source and target measures satisfy \cref{ass:source_target_measures}. Let $\gamma_{\eps,\lambda}$ be the minimizer of $K_{\eps,\lambda}(\mu,\nu)$. As $\lambda\to\infty$, $\gamma_{\eps,\lambda}\rightharpoonup \gamma_{\eps}$, where $\gamma_\eps$ is the optimal plan from $\mu$ to $\nu$ under the $c_\eps(x,y)$ cost function. Then as $\eps\to 0$, $\gamma_{\eps}\rightharpoonup \gamma_{\KR}$, the Knothe-Rosenblatt plan from $\mu$ to $\nu$.

\end{theorem}

\begin{proof}
    By \cref{lem:soft-hard}, the soft-constraint minimal value is bounded above by the hard-constraint minimal value. Dividing through by $\lambda$
    \[
 \mathcal{D}_{K L}\left((\pi_Y)_\sharp \gamma_{\eps,\lambda} \| \nu\right)+\frac{1}{\lambda} \int_{\mathbb{R}^d\times\R^d} c_\eps(x,y) d\gamma_{\eps,\lambda}
            \leq \frac{1}{\lambda}\int_{\mathbb{R}^d\times\R^d} c_\eps(x,y) d\gamma_{\eps} .\] Taking $\lambda\to \infty$, we get
\[
\lim_{\lambda\to \infty} \mathcal{D}_{K L}\left((\pi_Y)_\sharp \gamma_{\eps,\lambda} \| \nu\right)=0.
\]
 This means that $\nu_{\eps,\lambda}:=(\pi_Y)_\sharp \gamma_{\eps,\lambda}\rightharpoonup\nu$ as $\lambda\to \infty$. So $\nu_{\eps,\lambda}\rightharpoonup\nu$ as $\lambda\to\infty$. Also, by \cref{cor:unique min}, $\gamma_{\eps,\lambda}$ is the optimal transport plan from $\mu$ to $\nu_{\eps,\lambda}$ with cost $c_\eps(x,y)$.  Then we apply the stability of optimal transport plans to conclude that $\gamma_{\eps,\lambda}\rightharpoonup\gamma_\eps$ when $\lambda\to 0$. Finally \cref{thm0} tells us that sending $\eps\to 0$ gives $\gamma_\eps\rightharpoonup \gamma_{\KR}$.
\end{proof}

\section{Knothe-Rosenblatt Velocity Fields}\label{section3}

Our second set of results deals with the %
analog of Knothe-Rosenblatt maps 
for the dynamic optimal transport problem. Here, we seek \emph{triangular} velocity fields where the Jacobian of the velocity is an upper triangular matrix. In many modern applications, it is preferable to work with velocity fields instead of transport maps. For example, flow-based generative models seek a dynamical system whose flow map transports samples from a reference distribution to a target distribution \cite{lipman2023flow, albergo2023building}. However, computing an exact solution to the dynamical OT problem whose flow map samples exactly from the target distribution can be challenging and often intractable. In practice, such a velocity field can be identified through a soft-constrained dynamic optimal transport problem~\cite{onken2021otflow, wang2023efficient}. In this section, we prove that the minimizer of a Benamou-Brenier (hard-constrained) problem with a weighed Euclidean norm converges to a triangular velocity field. Before presenting our main result in Theorem~\ref{thm1}, we define the dynamical OT problem below with the weighted norm $\|x\|_\eps=c_\eps(x,x)^{1/2}$ given in~\eqref{eq:weighted_cost}.

\begin{definition}[Hard-constrained Benamou-Brenier]
Let $(CE)$ denote the set of distributions and velocity fields satisfying the continuity equation \eqref{CE}. Given absolutely continuous source and target measures $\mu$ and $\nu$, let  %
    \begin{align}\label{eq:bb-weighted}
    C_\eps(\mu,\nu):=\inf_{\rho,v} \left\{\int_0^1\int_{\R^d}\|v_t(x)\|_\eps^2\ \rho_t(x)\ dxdt:(\rho,v)\in(CE),\rho_0(x)=\mu(x),\rho_1(x)=\nu(x)\right\}.
\end{align}
\end{definition}
We define $X_\eps(t,x)$ to be the flow solving \eqref{eq:flowmap} with optimal velocity $v_\eps$ that solves~\eqref{eq:bb-weighted}. 
Our main contribution in Theorem~\ref{thm1} is to prove the dynamic analog of Theorem 2.1 in~\cite{carlier2008knothes}: we show that the sequence of weighted Benamou-Brenier problems in~\eqref{eq:bb-weighted} define velocity fields that converge to a triangular velocity fields with upper triangular Jacobian in space. Hence, the limiting flow induces a Knothe-Rosenblatt map. 
While parametric approaches to compute triangular velocity fields have been proposed in~\cite{kerrigan2024dynamic}, \cite{chemseddine2024conditionalwassersteindistancesapplications}, to our knowledge, this is the first guarantee for computing triangular velocity fields though a continuation method for dynamic optimal transport. Supporting results for the proof of the main result in this section can be found in Appendix~\ref{appendix: C}.

\begin{theorem}[Convergence to a triangular velocity field]\label{thm1} 
Let $\Omega,\Omega'$ be connected bounded open subsets of $\R^n$. Suppose that $\mu(dx)=f(x)dx$ and $\nu(dy)=g(y)dy$, with the respective densities $f, g$ being bounded from above and below. Furthermore, assume that $f\in C^k(\overline\Omega)$, $g\in C^k(\overline{\Omega'})$ with $\Omega,\Omega' \in C^{k+2}$ for $k\in\N$. Let $(\rho_\eps,v_\eps)$ be the minimizer to $C_\eps(\mu,\nu)$. Then, the minimizing velocity field $v_\eps(t,x)$ converges to  the velocity field induced by the Knothe-Rosenblatt map $T_{\KR}$ from $\mu$ to $\nu$ given by $$v_{\KR}(t,x) = T_{\KR}(X_{\KR}^{-1}(t,x))-X^{-1}_{\KR}(t,x),$$
pointwise a.e. on $\R^d$ and $t$, where the map $X_{\KR}(t) := t T_{\KR} + (1-t)x$ for all $t$.
\end{theorem}

\begin{proof} 
Let $T_\eps(x)$ be the optimal transport map pushing $\mu$ and $\nu$ with the $c_\eps(x,y)$ cost. The minimizing velocity field takes the form 
\[
v_\eps(t,x) = T_\eps(X^{-1}_\eps(t,x))-X^{-1}_\eps(t,x)
\]
where $X_\eps(t,x) = tT_\eps(x) +(1-t)x$ and $X^{-1}_\eps(t,x)$ is the inverse in space, which exists owing to the fact that $T_\eps(x)$ is a diffeomorphism, and the function $X_\eps(t,\cdot):\overline\Omega\to \tilde\Omega := t\overline\Omega+(1-t)\overline{\Omega'}$ is also diffeomorphic. Our goal is to show that 
\[
v_\eps(t,x) = T_\eps(X^{-1}_\eps(t,x))-X^{-1}_\eps(t,x) \to v_{\KR}(t,x)=T_{\KR}(X^{-1}_{\KR}(t,x))-X^{-1}_{\KR}(t,x) 
\]
where $v_{\KR}(t,x)$ is a triangular velocity field. Firstly, recall that since $T_\eps\to T_{\KR}$ in $L^2(\mu)$, it follows that $X_{\eps}(t,x) = tT_\eps(x)+(1-t)x \to X_{\KR}(t,x):=tT_{\KR}(x)+(1-t)x$ in $L^2(\mu)$. Also, since $\rho^\eps_t = X_{\eps}(t)_\sharp\mu$, Lemma \ref{lem:xt opt} tells us $X_\eps(t,x)$ is the optimal transport map from $\mu:=\rho_0$ to $\rho^\eps_t$. Next, for any $\varphi\in C_b(\overline\Omega)$, 
     \begin{align*}
    \int_{\overline\Omega} \varphi(y) d\rho_t^\eps &= \int_{\overline\Omega} \varphi(y) d\left( X_{\eps}(t,\cdot)_\sharp\mu\right)\\
    &=\int_{\overline\Omega} \varphi (X_\eps(t,x)) d\mu \to \int \varphi( X_{\KR}(t,x))d\mu = \int \varphi(y)d\rho_t
    \end{align*}
  up to a subsequence. This is because $X_\eps(t,x)\to X_{\KR}(t,x)$ pointwise up to a subsequence and so $\varphi(X_\eps(t,x))\to \varphi(X_{\KR}(t,x))$. Moreover,  $\varphi(X_\eps(t,x))$ is bounded for all $\eps$, since $\varphi$ is bounded. Hence, we can apply the Dominated convergence theorem and we have $\rho^\eps_t\rightharpoonup\rho_t$ for $\rho_t := X_{\KR}(t,\cdot)_\sharp\mu$, up to a subsequence. Notice that 
\[X_{\KR}(t,x) = \left(tT_1(x_1)+(1-t)x_1,tT_2(x_1,x_2)+(1-t)x_2,\ldots,t T_d(x_1,\ldots,x_d)+(1-t)x_d\right).\]
By Lemma \ref{lem:xt opt}, $T^1_t(x_1):=tT_1(x_1)+(1-t)x_1$ is an optimal map between $\mu^1$ and $(T_t^1)_\sharp\mu^1$, $T^2_t(x_2):= tT_2(x_1,x_2)+(1-t)x_2$ is an optimal map between the conditionals $\mu_{x_1}$ and $(T^2_t)_\sharp\mu_{x_1}$, and so on. This implies that $X_{\KR}(t,\cdot)$ is the Knothe-Rosenblatt map from $\mu$ to $\rho_t:= X_{\KR}(t,\cdot)_\sharp\mu$.\\ 

\noindent We now analyze the convergence of $X^{-1}_\eps(t,\cdot)$. 
  Recall that $X_\eps(t,\cdot)$ is a diffeomorphic sequence of functions, and converges, up to a subsequence, almost uniformly to $X_{\KR}(t,\cdot)$--which is continuous. We also have that $\rho_t^\eps = X_\eps(t,\cdot)_\sharp\mu$ and recall that $d\mu(x) = f(x)dx$. Then, the density of $\rho_t^\eps$ given by
  \[\rho_t^\eps(y) = f(X_\eps^{-1}(t,y))|\mathrm{det}\ DX^{-1}_\eps(t,y)|dy\] is bounded below by a constant depending on $\eps$, since $f$ is uniformly bounded by assumption and $|\mathrm{det}\ DX^{-1}_\eps(t,\cdot)|>c(\eps)\geq 0$. Given that the measures are bounded below by a constant, we can apply Lemma \ref{lem: inverse conv} which says $X_\eps^{-1}(t,\cdot)\to X^{-1}_{\KR}(t,\cdot)$ a.e. on $\tilde\Omega.$ Then, $X^{-1}_{\KR}(t,\cdot)$ is a triangular map from $\rho_t$ to $\mu$.\\ 
  
  \noindent We now show  
    $v_\eps(t,x) \to v_{\KR}(t,x)$ 
    pointwise a.e.\thinspace on $\R^d$.  %
    Indeed, 
    \begin{align*}
        &|v_\eps(t,x) -v_{\KR}(t,x)|\\
        &= \left|T_\eps(X^{-1}_\eps(t,x))-X^{-1}_\eps(t,x)-T_{\KR}(X^{-1}_{\KR}(t,x))+X^{-1}_{\KR}(t,x)\right|\\
        &\leq \left|T_\eps(X^{-1}_\eps(t,x))-T_{\KR}(X^{-1}_{\KR}(t,x))\right|+ \left|X^{-1}_{\KR}(t,x)-X^{-1}_\eps(t,x)\right|\\
        &\leq \left|T_\eps(X^{-1}_\eps(t,x))-T_{\KR}(X^{-1}_\eps(t,x))\right|+\left|T_{\KR}(X^{-1}_\eps(t,x))-T_{\KR}(X^{-1}_{\KR}(t,x))\right|+ \left|X^{-1}_{\KR}(t,x)-X^{-1}_\eps(t,x)\right|.
    \end{align*}
    For the first term, we note that
    \begin{align*}
            &\int_{\R^d}\left|T_\eps(X^{-1}_\eps(t,\cdot))-T_{\KR}(X^{-1}_\eps(t,\cdot))\right|^2d\rho_t^\eps=\int_{\R^d}\left|T_\eps-T_{\KR}\right|^2d\mu,
    \end{align*}
    which goes to 0 when $\eps\to 0$. Hence, the first term goes to zero pointwise a.e. up to a subsequence. 
    For the second term, by the Lipschitz-property of $T_{\KR}$ from Corollary \ref{cor: TK Lip}, we have
    \begin{align*}
        &\left|T_{\KR}(X^{-1}_\eps(t,x))-T_{\KR}(X^{-1}_{\KR}(t,x))\right| \leq \mathrm{Lip}(T_{\KR})\left|X^{-1}_\eps(t,x) -X^{-1}_{\KR}(t,x)\right|.
    \end{align*}
    Then, the right hand side and the third term go to 0 when $\eps\to 0$ pointwise a.e., up to a subsequence, by the result above. %
    Thus, $v_\eps$ converges to a triangular velocity field up to a subsequence. %
\end{proof}
We now define the soft-constrained variant of the Benamou-Brenier problem in~\eqref{eq:bb-weighted}.
\begin{definition}[Soft-constrained Benamou-Brenier]
 Let $(CE)$ denote the set of distributions and velocity fields satisfying the continuity equation \eqref{CE}. Given $\epsilon > 0$, $\lambda > 0$, we define
    \begin{align}\label{eq:mfg}
 C_{\eps,\lambda}(\mu,\nu)=\min _{\rho,v}&\bigg\{\lambda \mathcal{D}_{KL}\left(\rho(1,\cdot) \| \nu\right) + \int_0^1 \int_{\mathbb{R}^d}\|v(t,x)\|_\eps^2 \ \rho(t, x) \ d x d t:
 (\rho,v)\in(CE),\ \rho(0,x)=\mu(x)\bigg\}.
 \end{align}
 \end{definition}
\begin{rem}
    The main difficulty in adapting the proof of Theorem~\ref{thm1} to prove that the minimizing velocity of the soft-constrained Benamou-Brenier problem converges to a triangular velocity field is that the convergence of $X^{-1}_{\eps,\lambda}(t,x)$ to $X^{-1}_{\KR}(t,x)$ is not guaranteed. This is because we have no control over the regularity of the density of $\rho_{\eps,\lambda}(1,\cdot)$, the terminal condition of the soft-constrained problem. Therefore, we cannot assert that $X_{\eps,\lambda}(t,x)$ is diffeomorphic, so Lemma~\ref{lem: inverse conv}, which is used in the proof of Theorem~\ref{thm1} for the hard-constraint problem, does not apply. %
\end{rem}

\section{Discussion}
\label{sec:discussion}

This work opens several directions that can be used to %
construct triangular transport maps. From a computational perspective, it is valuable to investigate efficient computational algorithms and methodologies to approximate triangular maps using the soft-constrained optimization problems studied in this work. In particular, understanding the rates at which the approximation parameters $\epsilon$ in the cost, and $\lambda$ for the constraint, can be balanced against each other to achieve optimal trade-offs in approximation error remains an open problem, especially when the source and target distributions are only prescribed empirically by collections of samples. As a starting point, the trade-off can be studied when the source and target distributions are multivariate Gaussian distributions, for which the transport maps should admit closed-form solutions. A separate line of exploration is to evaluate different triangular maps that arise for each variable ordering. Considering alternative weightings in the cost function will change the orderings of variables in the limiting map and can be used to capture sparsity or simple functional forms in the map's components~\cite{spantini2018inference}. As a result, one can optimize the scalings in the cost function to effectively discover the variable orderings that most succinctly describe the underlying structure of the data.

\section{Acknowledgments}

RB acknowledges support from the von K\'{a}rm\'{a}n instructorship at Caltech,
the Air Force Office of Scientific Research MURI on “Machine Learning and Physics-Based Modeling
and Simulation” (award FA9550-20-1-0358) and a Department of Defense (DoD) Vannevar Bush
Faculty Fellowship (award N00014-22-1-2790) held by Andrew M. Stuart. MN acknowledges support from the Caltech SURF program that facilitated this research, Jun Kitagawa and Lauren Conger for many helpful discussions and user Vim from Mathematics Stack Exchange for the reference \cite{barvinek1991convergence}. BZ was supported by the AFOSR grant FA9550-21-1-0354. FH is supported by start-up funds at the California Institute of Technology and by NSF CAREER Award 2340762.

\appendix

\section{Proof of \cref{thm:kr stab}}

\label{app:main_proofs}

    Let $\gamma_\varepsilon$ denote the optimal coupling between $\mu_\varepsilon$ and $\nu_\varepsilon$ under cost $c_\varepsilon(x,y)$. Denote $\gamma_{\KR, \varepsilon}$ to be the Knothe-Rosenblatt coupling between $\mu_\varepsilon$ and $\nu_\varepsilon$. We closely follow the proof of Theorem 2.23 in \cite{otam}, but we highlight the main differences which rely on the application of the stability of optimal transport and KR plans under converging sequences of the source and target measures. 
    \begin{proof}
        We first show that $\gamma$, the limit of the optimal couplings $\gamma_\epsilon$ shares the same $(x_1,y_1)$ marginal as $\gamma_{\KR}$, the limit of the KR couplings $\gamma_{\KR, \epsilon}$. Then through disintegration with respect to $(x_1,y_1)$, we show that the conditionals of $(x_2,y_2)$ given $(x_1,y_1)$ of $\gamma$ and $\gamma_{\KR}$ are the same, which implies all the couplings between $(x_1,x_2,y_1,y_2)$ are equivalent. This argument is then repeated $d$ times to show all the marginals of $\gamma_{\varepsilon}$ and $\gamma_{\KR}$ coincide. In what follows, recall the notation for marginals and conditionals introduced in \cref{sec:notat}.  
   
\paragraph{Equivalence of the $(x_1,y_1)$ marginals}
By the optimality of $\gamma_\varepsilon$ we have the inequality
\begin{align}\label{eq:optimalityKRineq}
    \int_{\R^d \times \R^d} c_\epsilon(x,y) d\gamma_\epsilon  \le \int_{\R^d \times \R^d} c_\epsilon(x,y) d \gamma_{\KR, \epsilon}.
\end{align}
Note that $c_\eps(x,y)\to (x_1-y_1)^2$ locally uniformly. By stability of the optimal couplings (\cref{thm stab}) we have $\gamma_\epsilon \rightharpoonup \gamma$, where $\gamma$ is an optimal plan from $\mu$ to $\nu$ with respect to $c_0(x,y) = (x_1-y_1)^2$. Moreover, by narrow convergence up to a subsequence of $\gamma_{\KR,\epsilon} \rightharpoonup \gamma_{\KR}$ (\cref{lem:closed KR}), and the fact that $\gamma_\epsilon$ and $\gamma_{\KR,\eps}$ are uniformly compactly supported, taking $\eps\to 0 $ of \eqref{eq:optimalityKRineq} yields up to a subsequence, 
\begin{align}
    \int_{\R^d \times \R^d} (x_1-y_1)^2 d\gamma \le \int_{\R^d \times \R^d} (x_1-y_1)^2 d\gamma_{\KR}. 
\end{align}
Let $\gamma^1\coloneqq (\pi^{x,y}_1)_\sharp \gamma$ and $\gamma^1_{\KR}\coloneqq (\pi^{x,y}_1)_\sharp \gamma_{\KR}$, which are projections of the $\gamma$ and $\gamma_{\KR}$ couplings on the $(x_1,y_1)$ marginal, respectively. Since $\gamma_{\KR}^1$ is constructed to be optimal with respect to the quadratic cost $c_0(x,y) = (x_1-y_1)^2$, we have that 
\begin{align}
    \int_{\R\times \R} (x_1-y_1)^2 d\gamma_{\KR}^1 \le \int_{\R\times \R} (x_1 - y_1)^2 d\gamma^1 \le \int_{\R \times \R}(x_1 - y_1)^2 d\gamma_{\KR}^1,
\end{align}
which implies that $\gamma^1 = \gamma_{\KR}^1.$

\paragraph{Base case: equivalence of the $\gamma^{1:2}$ couplings}
 Let $\gamma^1_\epsilon = (\pi^{x,y}_1)_\sharp \gamma_\epsilon$ and $\gamma^1_{\KR,\epsilon} = (\pi^{x,y}_1)_\sharp \gamma_{\KR,\epsilon}$ respectively denote the projection of the optimal $\gamma_\epsilon$ and KR $\gamma_{\KR,\epsilon}$ couplings on the $(x_1,y_1)$ marginal. Recall that $\gamma_{\KR, \epsilon}$ and $\gamma_\epsilon$ have the same $x$ and $y$ marginals as they are both couplings between $\mu_\epsilon$ and $\nu_\epsilon$. Denote $\mu_\epsilon^1 = (\pi^x_1)_\sharp \mu_\epsilon$, $\nu_\epsilon^1 = (\pi^y_1)_\sharp \nu_\epsilon$ to be the $x_1$ and $y_1$ marginals of $\mu_\epsilon$ and $\nu_\epsilon$ respectively. Since the first component of the KR coupling $\gamma_{\KR, \varepsilon}^1$ is constructed to be the optimal between $\mu_\epsilon^1$ and $\nu_\epsilon^1$ with respect to the cost $c_0(x,y) = (x_1 - y_1)^2$, the projection of the optimal coupling on the $(x_1,y_1)$ marginals $\gamma^1_\varepsilon$ cannot be more optimal so we have the inequality
\begin{align}\label{eq:baseinequality}
    \int_{\R \times \R} (x_1 - y_1)^2 d\gamma_{\KR, \epsilon}^1 \le \int_{\R\times \R} (x_1 - y_1)^2 d \gamma_{\epsilon}^1.
\end{align}
Applying this inequality to \eqref{eq:optimalityKRineq}, we have
\begin{align*}
    &\int_{\R^d\times\R^d} (x_1-y_1)^2 d\gamma_{\KR,\eps} + \int_{\R^d\times\R^d}\left[\eps(x_2-y_2)^2 +\cdots + \eps^{d-1}(x_d-y_d)^2\right]d\gamma_{\eps}\\
    &\leq \int_{\R^d\times\R^d}c_\eps(x,y)d\gamma_{\eps} \leq \int_{\R^d\times\R^d}\left[ (x_1-y_1)^2 + \eps(x_2-y_2)^2 +\cdots + \eps^{d-1}(x_d-y_d)^2\right]d\gamma_{\KR,\eps}.
\end{align*}
Subtracting common terms, dividing by $\epsilon$, sending $\epsilon \to 0$, and appealing again to \cref{thm stab,lem:closed KR}, we obtain the inequality
\begin{align}\label{eq:marginal2ineq}
    \int_{\R^2 \times \R^2} (x_2-y_2)^2 d\gamma^{1:2} \le \int_{\R^2 \times \R^2} (x_2 - y_2)^2 d\gamma_{\KR}^{1:2}.
\end{align}
By disintegration and by the fact that $\gamma^1 = \gamma^1_{\KR}$, we may write the above inequality as follows
\begin{align}\label{eq:24 in FH notes}
    \int_{\R \times \R} d\gamma^1 \int_{\R\times\R} (x_2 - y_2)^2 d\gamma^{1:2}_{(x_1,y_1)} \le \int_{\R\times \R} d\gamma_{\KR}^{1} \int_{\R\times \R} (x_2 - y_2)^2 d\gamma_{\KR (x_1,y_1)}^{1:2}.
\end{align}
We want to check that the conditional measures $\gamma^{1:2}_{(x_1,y_1)}$ and $\gamma^{1:2}_{\KR (x_1,y_1)}$ have the same source and target marginals, meaning that we need to show that for any test function $\varphi$, 
\begin{align}
    \int_{\R\times \R} \varphi(x_2) d\gamma^{1:2}_{(x_1,y_1)}(x_2,y_2) = \int_{\R^\times\R} \varphi(x_2) d\gamma^{1:2}_{\KR (x_1,y_1)}(x_2,y_2).
\end{align}
To show this is true for $\gamma^1$-almost any $(x_1,y_1)$, it suffices to show that for any test functions $\varphi$ and $\psi$
\begin{align*}
    &\int_{\R\times \R} \psi(x_1,y_1) d\gamma^1(x_1,y_1) \int_{\R\times \R} \varphi(x_2) d\gamma^{1:2}_{(x_1,y_1)}(x_2,y_2) \\
    &=\int_{\R\times \R} \psi(x_1,y_1) d\gamma^1_{\KR}(x_1,y_1) \int_{\R\times \R} \varphi(x_2) d\gamma^{1:2}_{\KR (x_1,y_1)} (x_2,y_2)\\
    \iff & \int_{\R^2\times \R^2} \psi(x_1,y_1) \varphi(x_2) d\gamma^{1:2} = \int_{\R^2\times \R^2}  \psi(x_1,y_1)\varphi(x_2) d\gamma^{1:2}_{\KR}.
\end{align*}
Now, to show this is true, we first notice that while the integrands are functions of three variables, i.e., $\psi(x_1,y_1)\varphi(x_2)$, we may write it as a function of only $(x_1,x_2)$. Specifically, $\gamma^1$ corresponds to an optimal map $T_1:\R\to\R$ that is monotone and invertible where $y_1 = T_1(x_1)$. Therefore, we only need to show
\begin{align} \label{eq:integraleq}
    \int_{\R^2\times \R^2} \psi(x_1,T_1(x_1)) \varphi(x_2) d\gamma^{1:2} = \int_{\R^2\times \R^2}  \psi(x_1,T_1(x_1))\varphi(x_2) d\gamma^{1:2}_{\KR}.
\end{align}
Next, since $\gamma^{1:2}$ and $\gamma^{1:2}_{\KR}$ by construction share the same source and target marginals $\mu^{1:2} = (\pi^{x,y}_{1:2})_\sharp \mu$ and  $\nu^{1:2} = (\pi^{x,y}_{1:2})_\sharp \nu$, expectations of functions of $(x_1,x_2)$ are immediately equal, that is \eqref{eq:integraleq} holds. 

A similar argument can be made for test function $\varphi$ that have arguments $y_2$ where we instead use the change of variables $x_1 = T_1^{-1}(y_1)$, which is well-defined since $T_1$ is a monotone invertible function. This allows us to conclude that $\gamma^{1:2}_{(x_1,y_1)}$ and $\gamma^{1:2}_{\KR(x_1,y_1)}$ are $\gamma^1$-a.e. equal. Moreover, since the conditional distribution $\gamma^{1:2}_{\KR(x_1,y_1)}$ is constructed to be optimal with respect to the cost function $(x_2 - y_2)^2$ and using \eqref{eq:24 in FH notes}, we have
\begin{align}
   \int_{\R\times \R} d\gamma_{\KR}^1 \int_{\R\times \R} (x_2 - y_2)^2 d\gamma_{\KR (x_1,y_1)}^{1:2} & \le \int_{\R \times \R} d\gamma^1 \int_{\R\times\R} (x_2 - y_2)^2 d\gamma^{1:2}_{(x_1,y_1)} \\ & \le \int_{\R\times \R} d\gamma_{\KR}^1 \int_{\R\times \R} (x_2 - y_2)^2 d\gamma_{\KR (x_1,y_1)}^{1:2}, \nonumber 
\end{align}
which implies $\gamma^{1:2} = \gamma_{\KR}^{1:2}$. 

\paragraph{Inductive step: Constructing a sub-optimal coupling }
For the inductive step, suppose we have $\gamma^{1:n} = \gamma_{\KR}^{1:n}$. We will prove $\gamma^{1:n+1} = \gamma^{1:n+1}_{\KR}$ and we do so similarly to the base case by first establishing a similar inequality to \eqref{eq:marginal2ineq}. However we cannot proceed as we did for the base case through an inequality such as \eqref{eq:baseinequality} as $\gamma^{1:n}_{\KR, \epsilon}$ is not optimal for cost $(x_1-y_1)^2 + \cdots + \epsilon^{n-1}(x_{n}-y_{n})^2$. Instead, we will establish the inequality by construction of a suboptimal measure.

Consider a suboptimal measure that consists of two components. The first is the optimal coupling $\eta^{1:n}_{ \epsilon} \in \mathcal{P}(\R^n \times\R^n)$ with source and target measures $\mu_\epsilon^{1:n}$ and $\nu_\epsilon^{1:n}$ that is optimal with respect to cost $(x_1 - y_1)^2 + \cdots + \epsilon^{n-1}(x_n-y_n)^2$. For the second component, given a fixed $(x_{1:n},y_{1:n})$, we consider the family of KR couplings $\zeta_{\KR,\eps (x_{1:n},y_{1:n})} \in \mathcal{P}(\R^{d-n}\times \R^{d-n})$ between the conditional measures $\mu_{\epsilon, x_{1:n}}$ and $\nu_{\epsilon, y_{1:n}}$. The conditional measures $\mu_{\eps,x_{1:n}}$ and $\nu_{\eps,y_{1:n}}$ converge narrowly to $\mu_{x_{1:n}}$ and $\nu_{y_{1:n}}$ by \cref{prop:cv of disintegration}. The suboptimal measure $\lambda_{\epsilon,n} \in \mathcal{P}(\R^d \times \R^d)$ is then defined as
\begin{align}
    \lambda_{\epsilon,n} = \eta^{1:n}_\epsilon \otimes \zeta_{\KR,\eps (x_{1:n},y_{1:n})}.
\end{align}
Here, $\lambda_{\epsilon,n}$ can be understood as a joint probability distribution that is constructed as a product of a marginal distribution $\eta_\epsilon^{1:n}$ and a conditional distribution $\zeta_{\KR(x_{1:n},y_{1:n})}.$ We first check that $\lambda_{\epsilon,n}$ is indeed a coupling between $\mu_\epsilon$ and $\nu_\epsilon$. Let $\varphi$ be an arbitrary Borel measurable bounded function on $\R^d$ and observe that
\begin{align*}
    \int_{\R^d\times \R^d} \varphi(x) d\lambda_{\epsilon,n} & = \int_{\R^n \times \R^n} \int_{\R^{d-n}\times \R^{d-n}} \varphi(x) d\zeta_{\KR,\eps(x_{1:n},y_{1:n}) }  d\eta_\epsilon^{1:n}  \\
    & = \int_{\R^n} \int_{\R^{d-n}} \varphi(x) d\mu_{\epsilon, x_{1:n}} d \mu_\epsilon^{1:n}  = \int_{\R^n} \varphi(x) d\mu_\epsilon.
\end{align*}
Here, the last equality is by the disintegration theorem. A similar argument can be made for functions of $y \in \R^d$, hence showing $\lambda_{\epsilon,n}$ has source and target measures $\mu_\epsilon$ and $\nu_\epsilon$. 
Also, since the conditional measure of $\mu_{\eps,x_{1:n}}$ and $\nu_{\eps,y_{1:n}}$ converge by \cref{prop:cv of disintegration}, we can apply \cref{lem:closed KR} to say that $$\zeta_{\KR,\eps(x_{1:n},y_{1:n})}\rightharpoonup\zeta_{\KR,(x_{1:n},y_{1:n})}.$$ This is the KR plan between $\mu_{x_{1:n}}$ and $\nu_{y_{1:n}}$. Next, by \cref{thm stab}, we know that $\eta^{1:n}_{\epsilon} \rightharpoonup  \gamma^{1:n}$, which is $\gamma_{\KR}^{1:n}$ by the induction hypothesis. Then, we can say that as we send $\eps\to 0$, we obtain convergence of the suboptimal coupling
\begin{equation}\label{eq:star}
\lambda_{\eps,n}=\eta^{1:n}_\eps\otimes \zeta_{\KR,\eps(x_{1:n},y_{1:n})}\rightharpoonup\gamma_{\KR}^{1:n}\otimes \zeta_{\KR,(x_{1:n},y_{1:n})} = \gamma_{\KR} 
\end{equation}
\\
\textbf{Inequality for the $(n+1)$st marginal.}
Now the proof proceeds similarly to the base step. 
First consider the following inequality analogous to \eqref{eq:marginal2ineq}. The first $n$ joint marginals of the suboptimal coupling $\eta_{\epsilon}^{1:n} = (\pi^{x,y}_{1:n})_\sharp\lambda_{\epsilon,n}$ is optimal between $\mu_\epsilon^{1:n}$ and $\nu_\epsilon^{1:n}$ with respect to cost $\sum_{i = 1}^n \epsilon^{i-1} (x_i - y_i)^2$. Therefore, the projection of the optimal coupling $\gamma^\epsilon$ on the $(x_{1:n},y_{1:n})$ marginals, $\gamma_\epsilon^{1:n}$,  cannot be more optimal, so we have the inequality
\begin{align} \label{eq:inductivemarginalinequality}
  \int_{\R^d \times \R^d } &\sum_{i = 1}^n \epsilon^{i-1} (x_i - y_i)^2 d\lambda_{\epsilon,n} =   \int_{\R^n \times \R^n} \sum_{i = 1}^n \epsilon^{i-1}(x_i - y_i)^2 d \eta_{\epsilon,n}  \\ \le &\int_{\R^n \times \R^n} \sum_{i =1 }^n \epsilon^{i-1} (x_i - y_i)^2 d \gamma^{1:n}_\epsilon  = \int_{\R^d \times \R^d } \sum_{i = 1}^n \epsilon^{i-1} (x_i - y_i)^2 d\gamma_\epsilon. \nonumber
\end{align}
By optimality of the optimal coupling $\gamma_\epsilon$ with respect to $c_\epsilon(x,y)$, we have
\begin{align*}
    \int_{\R^d\times \R^d} c_\epsilon(x,y) d\gamma_\epsilon \le \int_{\R^d\times \R^d} c_\epsilon(x,y) d\lambda_{\epsilon,n}.
\end{align*}
Combining this with \eqref{eq:inductivemarginalinequality}, we have
\begin{align*}
   \int_{\R^d \times \R^d} & \sum_{i = 1}^n \epsilon^{i-1} (x_i - y_i)^2 d\lambda_{\epsilon,n} +  \int_{\R^d\times \R^d} \sum_{i = n+1}^d \epsilon^{i - 1} (x_i - y_i)^2 d\gamma_\epsilon \\ &  \le \int_{\R^d \times \R^d} c_\epsilon(x,y) d\gamma_\epsilon\le \int_{\R^d \times \R^d} \sum_{i = 1}^n \epsilon^{i -1} (x_i - y_i)^2 d\lambda_{\epsilon,n} + \int_{\R^d \times \R^d} \sum_{i = n+1}^d \epsilon^{i -1} (x_i - y_i)^2 d\lambda_{\epsilon,n}. 
\end{align*}
Cancelling common terms, dividing by $\epsilon^{n}$, sending $\epsilon \to 0$, recalling the limit~\eqref{eq:star} and that $\gamma_\eps^{1:n+1}\rightharpoonup \gamma^{1:n+1}$, and marginalizing out all variables above index $n+1$, we have
\begin{align} \label{eq:np1marginalienquality}
    \int_{\R^{n+1} \times \R^{n+1}}(x_{n+1} - y_{n+1})^2 d\gamma^{1:n+1} &\le  \int_{\R^{n+1} \times \R^{n+1}}(x_{n+1} - y_{n+1})^2 d(\gamma_{\KR}^{1:n} \otimes \zeta_{\KR(x_{1:n},y_{1:n})}^{n+1}).
\end{align}
Here, recall that $\zeta_{\KR (x_{1:n},y_{1:n})}^{n+1}$ is the  $(x_{n+1}, y_{n+1})$ conditional-marginal of given $(x_1,\ldots, x_n, y_1, \ldots, y_n)$. 

By disintegration and by the inductive hypothesis that $\gamma^{1:n} = \gamma_{\KR}^{1:n}$,
implies that 
\begin{align} \label{eq:np1disintegration}
    \int_{\R^n \times \R^n } d\gamma^{1:n} \int_{\R \times \R} (x_{n+1} - y_{n+1})^2 d \gamma^{1:n+1}_{(x_{1:n}, y_{1:n})} \le \int_{\R^n \times \R^n} d\gamma^{1:n}_{\KR} \int_{\R \times \R} (x_{n+1} - y_{n+1})^2 d \zeta_{\KR(x_{1:n},y_{1:n})}^{n+1}.
\end{align}
Now we check that couplings $\gamma_{(x_{1:n},y_{1:n})}^{1:n+1}$ and $\zeta_{\KR(x_{1:n},y_{1:n})}^{n+1}$ have the same source and target marginals, meaning that we show for any test function $\varphi$ of $x_{n+1} \in \R$,
\begin{align}
    \int_{\R\times \R} \varphi(x_{n+1}) d\gamma^{1:n+1}_{(x_{1:n}, y_{1:n})} = \int_{\R\times \R} \varphi(x_{n+1}) d\zeta^{n+1}_{\KR (x_{1:n}, y_{1:n})}
\end{align}
for $\gamma^{1:n}$-a.e. $(x_{1:n}, y_{1:n})$. Recall that $\zeta_{\KR(x_{1:n},y_{1:n})}^{n+1}$ is the $(n+1)$st conditional-marginal of $\gamma^{1:n+1}_{\KR (x_{1:n}, y_{1:n})}$, ie.e $\gamma_{\KR}^{1:n+1} = \gamma_\KR^{1:n}\otimes \zeta^{n+1}_{\KR(x_{1:n},y_{1:n})}$. Hence it suffices to show that for every function $\psi$ of $(x_{1:n}, y_{1:n})$, 
\begin{align}
    \int_{\R^{n+1}\times \R^{n+1}}\psi(x_{1:n},y_{1:n}) \varphi(x_{n+1}) d\gamma^{1:n+1} =     \int_{\R^{n+1}\times \R^{n+1}}\psi(x_{1:n},y_{1:n}) \varphi(x_{n+1}) d\gamma^{1:n+1}_{\KR }.
\end{align}
 Invoking the inductive hypothesis again, the optimal coupling (which is a KR coupling) $\gamma^{1:n} = \gamma_{\KR}^{1:n}$ induces a deterministic invertible optimal map $y_{1:n} = T_{1:n}(x_{1:n})$, where $$T_{1:n} = \begin{bmatrix} T_1(x_1), T_2(x_1,x_2), \ldots, T_n(x_{1:n})\end{bmatrix}^\top $$ and component $T_i$ is monotone and invertible with respect to $x_i$ given any fixed $x_{1:i-1}$. This implies that we only need to show that
\begin{align}
    \int_{\R^{n+1}\times \R^{n+1}}\psi(x_{1:n},T_{1:n}(x_{1:n})) \varphi(x_{n+1}) d\gamma^{1:n+1} =     \int_{\R^{n+1}\times \R^{n+1}}\psi(x_{1:n},T_{1:n}(x_{1:n})) \varphi(x_{n+1}) d\gamma^{1:n+1}_{\KR },
\end{align}
which we know to be true since $\gamma^{1:n+1}$ and $\gamma_{\KR}^{1:n+1}$ has the same source and target measures $\mu^{1:n+1}$ and $\nu^{1:n+1}$ by the hypothesis of the theorem. A similar argument can be made for functions $\varphi(y_{n+1})$, in which we use the fact that the function $T_{1:n}$ has a well-defined inverse function. 
Returning to \eqref{eq:np1disintegration}, we use the fact that for a given $(x_{1:n},y_{1:n})$, $\zeta_{\KR(x_{1:n},y_{1:n})}^{n+1}$ is the optimal coupling between conditional-marginals $\mu_{x_{1:n}}^{1:n+1}$ and $\nu_{y_{1:n}}^{1:n+1}$, to obtain
\begin{align}
        \int_{\R^n \times \R^n} d\gamma^{1:n}_{\KR} \int_{\R \times \R} (x_{n+1} - y_{n+1})^2 d \zeta_{\KR(x_{1:n},y_{1:n})}^{n+1} &\le \int_{\R^n \times \R^n } d\gamma^{1:n} \int_{\R \times \R} (x_{n+1} - y_{n+1})^2 d \gamma^{1:n+1}_{(x_{1:n}, y_{1:n})} \\ & \le \int_{\R^n \times \R^n} d\gamma^{1:n}_{\KR} \int_{\R \times \R} (x_{n+1} - y_{n+1})^2 d \zeta_{\KR(x_{1:n},y_{1:n})}^{n+1}\nonumber,
\end{align}
which allows us to conclude $\gamma^{1:n+1} = \gamma^{1:n+1}_\KR$. This argument is then iterated until $n+1 = d$. 
\end{proof}

\section{Supporting results for Section 3} \label{app:additional_proofs}

\begin{lemma}[Lemma 4.4 in \cite{Villani2008OptimalTO}]\label{lem:vil}
   Let $P,Q$ be tight subsets of $\mathcal{P}(\R^d)$. Then the set $\Pi(P, Q)$ of all couplings whose marginals lie in $P$ and $Q$ respectively, is itself tight in $\mathcal{P}(\R^d \times \R^d)$. 
\end{lemma}

\begin{proposition}\label{prop: exists}
    There exists a minimizer for $K_{\eps,\lambda}(\mu,\nu)$. 
\end{proposition}
\begin{proof} 
Pick out a minimizing sequence $(\gamma_i)_{i\in \N}$ of $K_{\eps,\lambda}(\mu,\nu)$. Recall that $\mu$ has finite second moment. Also, since for all $i\in \N$, $(\pi_Y)_\sharp\gamma_i \ll \nu$ and $\nu$ has compact support, this means that $(\pi_Y)_\sharp\gamma_i$ has uniform finite second moment. Thus $\mathcal{Q}=((\pi_Y)_\sharp\gamma_i)_{i\in\N}$ has uniformly bounded second moment. So it is tight. Then \cref{lem:vil} tells us that $\Pi(\mu,\mathcal{Q})$ is tight, which means that there exists a cluster point for the minimizing sequence $(\gamma_i)_{i\in\N}$. This cluster point is the minimizer for $K_{\eps,\lambda}(\mu,\nu)$.
\end{proof}

\begin{proposition}\label{cor:unique min}
Let $\mu$ be a probability measure with finite second moment such that $\mu\ll dx$, $\nu$ be a probability measure with finite second moment. Let $\gamma_{\eps,\lambda}$ be a minimizer of the soft-constrained problem $K_{\eps,\lambda}(\mu,\nu)$. Then $\gamma_{\eps,
\lambda}$ is the minimizer of the hard-constrained problem $K_\eps(\mu,\nu_{\eps,\lambda})$ where $\nu_{\eps,\lambda}:= (\pi_Y)_\sharp\gamma_{\eps,\lambda}.$ Furthermore, $\gamma_{\eps,\lambda}$ is unique, and $\gamma_{\eps,\lambda} = (id,T_{\eps,\lambda})_\sharp\mu$, where $T_{\eps,\lambda}$ is the optimal transport map between $\mu$ and $\nu_{\eps,\lambda}.$
\end{proposition}
\begin{proof}
 Let $\gamma_{\eps,\lambda}$ be a minimizer of $K_{\eps,\lambda}(\mu,\nu)$, $\nu_{\eps,\lambda}:=(\pi_Y)_\sharp\gamma_{\eps,\lambda}$ and 
    \begin{equation}\label{eq:admis plan}
    \zeta_{\eps,\lambda}:=\argmin_{\gamma}\left\{\int_{\R^d\times\R^d} c_\eps(x,y)d\gamma:\gamma\in\Pi\left(\mu,\nu_{\eps,\lambda}\right)\right\}.
    \end{equation}
   We would like to show that $\zeta_{\eps,\lambda} = \gamma_{\eps,\lambda}$. Suppose for contradiction that $\zeta_{\eps,\lambda}\neq \gamma_{\eps,\lambda}$. Then
    \[
    \int_{\R^d\times\R^d} c_\eps(x,y)d\zeta_{\eps,\lambda}<\int_{\R^d\times\R^d} c_\eps(x,y)d\gamma_{\eps,\lambda}
    \]
   as $\gamma_{\eps,\lambda}$ is an admissible plan for \eqref{eq:admis plan}. Since $\gamma_{\eps,
   \lambda}, \zeta_{\eps,\lambda}\in \Pi(\mu,\nu_{\eps,\lambda})$, we have 
    \[
    \mathcal{D}_{K L}\left((\pi_Y)_\sharp \gamma_{\eps,\lambda}\|\nu\right) = \mathcal{D}_{KL}\left((\pi_Y)_\sharp \zeta_{\eps,\lambda}\|\nu\right)=\mathcal{D}_{K L}\left(\nu_{\eps,\lambda}\|\nu\right).
    \]
    This implies that $\gamma_{\eps,\lambda}$ cannot be the minimizer for $K_{\eps,\lambda}(\mu,\nu)$ since
    \begin{align*}
     \lambda\mathcal{D}_{KL}\left( (\pi_Y)_\sharp \zeta_{\eps,\lambda}\|\nu\right) +\int_{\R^d\times\R^d} c_\eps(x,y)d\zeta_{\eps,\lambda}  <  \lambda\mathcal{D}_{K L}\left( (\pi_Y)_\sharp \gamma_{\eps,\lambda}\|\nu\right) +\int_{\R^d\times\R^d} c_\eps(x,y)d\gamma_{\eps,\lambda},
    \end{align*}
    a contradiction. Hence $\gamma_{\eps,\lambda}$ is the optimal transport plan from $\mu$ to $\nu_{\eps,\lambda}$. By Brenier's theorem, $\gamma_{\eps,\lambda}$ is unique, and $\gamma_{\eps,\lambda}=(id\times T_{\eps,\lambda})_\sharp\mu$, where $T_{\eps,\lambda}$ is the optimal map between $\mu$ and $\nu_{\eps,\lambda}$.
\end{proof}

\begin{proposition}\label{lem:soft-hard}
    The minimum of the soft-constraint problem is always bounded above by the minimum of the hard-constraint problem, that is $K_{\eps,\lambda}(\mu,\nu)\leq K_{\eps}(\mu,\nu)$.
\end{proposition}
\begin{proof}
    Let $\gamma_{\eps,\lambda}$ be a minimizer to the soft-constraint problem $K_{\eps,\lambda}(\mu,\nu)$ which exists and is unique thanks to \cref{prop: exists,cor:unique min}. Thus, for all couplings $\zeta$ with first marginal $\mu$, we get
    \begin{align*}
        \lambda\mathcal{D}_{K L}\left((\pi_Y)_\sharp \gamma_{\eps,\lambda} \| \nu\right)&+ \int_{\mathbb{R}^d\times\R^d} c_\eps(x,y) d\gamma_{\eps,\lambda}\\
           &\leq\lambda\mathcal{D}_{K L}\left((\pi_Y)_\sharp \zeta \| \nu\right)+\int_{\mathbb{R}^d\times\R^d} c_\eps(x,y) d\zeta.
    \end{align*}
    Take $\zeta=\gamma_\eps$, the minimizer to $\min_\gamma\{\int_{\R^d}c_\eps(x,y)d\gamma:\gamma\in\Pi(\mu,\nu)\}$. Since $\gamma_\eps\in \Pi(\mu,\nu)$,
    \begin{align*}
            \lambda\mathcal{D}_{K L}\left((\pi_Y)_\sharp \gamma_{\eps,\lambda} \| \nu\right)&+ \int_{\mathbb{R}^d\times\R^d} c_\eps(x,y) d\gamma_{\eps,\lambda}\\
            &\leq \lambda\mathcal{D}_{K L}\left((\pi_Y)_\sharp \gamma_{\eps} \| \nu\right)+ \int_{\mathbb{R}^d\times\R^d} c_\eps(x,y) d\gamma_{\eps}\\
            &=\int_{\mathbb{R}^d\times\R^d} c_\eps(x,y) d\gamma_{\eps}= \min_\gamma\{\int_{\R^d}c_\eps(x,y)d\gamma:\gamma\in\Pi(\mu,\nu)\}.
    \end{align*}
    This shows that $K_{\eps,\lambda}(\mu,\nu)\leq K_{\eps}(\mu,\nu)$.
\end{proof}

\begin{theorem}[Theorem 2.6.9 in \cite{figalli2021invitation}, Stability of Optimal Transport Plans]\label{thm stab}
    Let $c_k(x,y) \in C^0(\R^d \times \R^d)$ be a sequence of costs bounded from below, and let $\mu_k$ and $\nu_k$ be two sequences of probability measures with $\mu_k \rightharpoonup \mu$ and $\nu_k \rightharpoonup \nu$ such that $\mu$ and $\nu$ are probability measures. For each $k \in \mathbb{N}$, let $\gamma_k \in \Pi\left(\mu_k, \nu_k\right)$ be an optimal transport plan %
    with respect to the cost $c_k$, and assume that $c_k \rightarrow c$ locally uniformly. Then, up to a subsequence, the sequence $\gamma_k$ converges narrowly to an optimal transport plan $\gamma \in \Pi(\mu, \nu)$ with respect to the cost $c$.
\end{theorem}

Next, we present a key theorem from \cite{goggin1994convergence}; then we interpret the result in our context. 

\begin{theorem}[Theorem 2.1 in \cite{goggin1994convergence}, Condition for convergence of conditionals]\label{thm:goggin}
Let $S_1$ and $S_2$ be complete, separable metric spaces and let $P_n\in \mathcal{P}(\Omega_n)$ and $P\in \mathcal{P}(\Omega)$ with associated measurable functions (random variables) \[(X_n,Y_n):\Omega_n\to S_1\times S_2,\qquad(X,Y):\Omega \to S_1\times S_2\] such that $(X_n,Y_n)_\sharp P_n\rightharpoonup (X,Y)_\sharp P$. Note that the measurable functions are defined such that \[X_n:\Omega_n\to S_1,\quad Y_n:\Omega_n\to S_2,\quad X:\Omega\to S_1,\quad \text{and}\quad Y:\Omega\to S_2.\]  Assume that there exists a reference probability measure $Q_n$ such that $P_n\ll Q_n$ on $\sigma((X_n,Y_n))$ with Radon-Nikodym derivative $L_n\circ (X_n,Y_n)$. Further assume that $(X_n,Y_n)_\sharp Q_n = (X_n)_\sharp Q_n\otimes(Y_n)_\sharp Q_n$ and that $\bar\mu_n:=(\pi_X)_\sharp Q_n$ and $\bar\nu_n:=(\pi_Y)_\sharp Q_n$ such that $\bar\mu_n\times\bar\nu_n \rightharpoonup\bar\mu\times\bar\nu$, some other product measure. Define $Q:=(X,Y)_\sharp (\bar\mu\times\bar\nu)$ and assume that $(X_n,Y_n,L_n\circ(X_n,Y_n))_\sharp Q_n\rightharpoonup (X,Y,L\circ(X,Y))_\sharp Q$ where $L\circ(X,Y)$ is some function that is a probability density under the measure $Q$. Then 
\begin{itemize}
    \item $P\ll Q$ on $\sigma((X,Y))$, $dP/dQ = L\circ(X,Y)$
    \item and for $P_n = (\pi_Y)_\sharp P_n\otimes P_n^Y$ and $P =  (\pi_Y)_\sharp P\otimes P^Y$, then $P^Y_n\rightharpoonup P^Y$,
\end{itemize}
where $P^Y_n$ and $P^Y$ are conditional measures given $Y$.
\end{theorem}

In what follows, we use the same notation for marginals and conditionals as introduced at the end of \cref{sec:notat}.
\begin{theorem}[Stability of Conditionals]\label{prop:cv of disintegration}
    Let $\mu,\mu_n\in \mathcal{P}(\Lambda)$, where $\Lambda$ is a compact set in $\R^d$, and $\mu_n,\mu\ll dx$. 
    If $\mu_n\rightharpoonup\mu$ and the density $L_n$ of $\mu_n$ converges pointwise a.e. to the density $L$ of $\mu$, then $\mu_{n,x_{1:i}}\rightharpoonup \mu_{x_{1:i}}$. 
\end{theorem}

\begin{proof}
Let $\mu_n^{1:i}=(\pi_{1:i})_\sharp\mu_n$ and $\mu_n^{i+1:d}=(\pi_{i+1:d})_\sharp\mu_n$. We now check the necessary conditions to apply \cref{thm:goggin}. Using the same notation as in \cref{thm:goggin}, let $\Omega_n= \Omega =:\Lambda$, $S_1 := \pi_{i+1:d}( \Lambda)$, $S_2 := \pi_{1:i}(\Lambda)$ and take $X_n=X = \pi_{i+1:d}$ and $Y_n=Y = \pi_{1:i}$ to be the respective projections, define $P_n := \mu_n$, $P:=\mu$ and define $Q_n := dx_K$ as the normalized Lebesgue measure on a compact hypercube $K$ containing $\Lambda$. Finally, let $\bar\mu_n=\bar\mu:= dx_K^{i+1:d}$ and $\bar\nu_n=\bar\nu:=dx_K^{1:i}$. Now we verify each condition one by one
\begin{enumerate}
    \item $(X_n,Y_n)_\sharp P_n \rightharpoonup (X,Y)_\sharp P$: This follows immediately from the fact that the composition of any continuous bounded function with any permutation is again a continuous bounded function.
    \item $P_n \ll Q_n$: This condition is trivial since $P_n = \mu_n$, which is assumed to be absolutely continuous with respect to $dx$.
    \item $(X_n,Y_n)_\sharp Q_n = (X_n)_\sharp Q_n\otimes (Y_n)_\sharp Q_n$: Since $dx$ is independent, this condition follows.
    \item $\bar\mu_n\otimes \bar\nu_n \rightharpoonup \bar\mu\otimes\bar\nu$: Since $\bar\mu_n = \bar\mu$ and $\bar\nu_n = \bar\nu$, this condition follows.
    \item  Recall that $L_n,L$ are the densities of $\mu_n,\mu$. We want to show that $(X_n,Y_n,L_n\circ(X_n,Y_n))_\sharp Q_n\rightharpoonup (X,Y,L\circ(X,Y))_\sharp Q$ where $L\circ(X,Y)$ is some function that it is a probability density under the measure $Q  = dx_K$. This is equivalent to showing that for all $\varphi \in C_b(K\times [0,\infty))$,
\begin{equation}\label{eq: Ln condition}
    \int_{K}\varphi(x_{i+1:d},x_{1:i},L_n(x_{i+1:d},x_{1:i}))dx \to \int_{K}\varphi( x_{i+1:d},x_{1:i},L(x_{i+1:d},x_{1:i}))dx\,.
\end{equation}
Since $L_n\to L$ pointwise a.e. and $\varphi$ is continuous, therefore $\varphi(x_{i+1:d},x_{1:i},L_n(x_{i+1:d},x_{1:i}))\to \varphi(x_{i+1:d},x_{1:i},L(x_{i+1:d},x_{1:i}))$ a.e. We also have that $\varphi$ is bounded and in $L^1(K\times [0,\infty))$. Therefore we can use the dominated convergence theorem to conclude \eqref{eq: Ln condition}.
\end{enumerate}
 Now we can apply \cref{thm:goggin} to conclude $\mu_{n,x_{1:i}}\rightharpoonup \mu_{x_{1:i}}$.
\end{proof}

\begin{rem}\label{rmk:assAvsB}
\cref{prop:cv of disintegration} still applies in the case where $\mu_n, \mu$ are not compactly supported. Indeed, if our measures are absolutely continuous with respect to a Gaussian measure, then the Gaussian itself satisfies the independence property across coordinates and is a probability measure, so the assumptions of the theorem are still satisfied. As a consequence, we can replace \cref{ass:source_target_measures} with \cref{asm:ac gaussian} in the results making use of \cref{prop:cv of disintegration}.
\end{rem}

\begin{proposition}[Stability of KR plans]\label{lem:closed KR}
    Let $\mu,\nu,\mu_n,\nu_n\in \mathcal{P}(\Lambda)$ for some compact subset $\Lambda$ in $\R^d$, and let $\gamma_{\KR,n}$ be the Knothe-Rosenblatt plan from $\mu_n$ to $\nu_n$. Assume $\mu_n\rightharpoonup\mu$, $\nu_n\rightharpoonup\nu$ with $\mu_n\ll dx$ and $\nu_n\ll dy$. Further, assume the densities of $\mu_n$ converge pointwise a.e. to the density of $\mu$ (resp.$\nu_n$ and $\nu$), then up to a subsequence, $\gamma_{\KR,n}\rightharpoonup \gamma_{\KR}$, a Knothe-Rosenblatt plan from $\mu$ to $\nu$.
\end{proposition}

\begin{proof}
 Let $\gamma_{\KR,n}$ be the KR plan from $\mu_n$ to $\nu_n$ for the Euclidean cost. Then $\gamma_{\KR,n}$ is also an optimal plan with respect to the cost $(x_1-y_1)^2$. By the stability of optimal transport in \cref{thm stab}, $\gamma_{\KR,n}\rightharpoonup\gamma$, which is the optimal plan from $\mu$ to $\nu$ with respect to the cost $(x_1-y_1)^2$. We now want to show that $\gamma = \gamma_\KR$, the Knothe-Rosenblatt plan from $\mu$ to $\nu$. By the disintegration theorem applied to each measure $\mu,\nu,\mu_n,\nu_n$ after conditioning on the first $i\in\{2,\ldots,d-1\}$ variables, we have
 \begin{align*}
    \mu_n(x) &= \mu_n^{1:i}(x_{1:i}) \otimes \mu_{n,x_{1:i}}(x_{i+1:d}) \\
    \mu(x) &= \mu^{1:i}(x_{1:i}) \otimes \mu_{x_{1:i}}(x_{i+1:d}) \\
    \nu_n(y) &= \nu_n^{1:i}(y_{1:i}) \otimes \nu_{n,y_{1:i}}(y_{i+1:d}) \\
    \nu(y) &= \nu^{1:i}(y_{1:i}) \otimes \nu_{y_{1:i}}(y_{i+1:d}) 
 \end{align*}
    Define 
      \[
    \begin{dcases}
    \mu_{n,x_{1:i}}^{i+1}(x_{i+1}):=(\pi_{i+1})_\sharp\mu_{n,x_{1:i}}(x_{i+1:d})\\
    \mu^{i+1}_{x_{1:i}}(x_{i+1}):=(\pi_{i+1})_\sharp\mu_{x_{1:i}}(x_{i+1:d}).
    \end{dcases}\]
    and
    \[
    \begin{dcases}
    \nu_{n,y_{1:i}}^{i+1}(y_{i+1}):=(\pi_{i+1})_\sharp\mu_{n,y_{1:i}}(y_{i+1:d})\\
    \nu^{i+1}_{y_{1:i}}(y_{i+1}):=(\pi_{i+1})_\sharp\mu_{y_{1:i}}(y_{i+1:d}).
    \end{dcases}\]
    For fixed $(x_1,\ldots,x_i)$ and $(y_1,\ldots,y_i)$, $\gamma_{\KR,n}$ is the optimal plan from $\mu_{n,x_{i:i}}^{i+1}$ to $\nu_{n,y_{i:i}}^{i+1}$ with quadratic cost $(x_{i+1}-y_{i+1})^2$. By \cref{prop:cv of disintegration}, since $\mu_n,\mu\ll dx$ and $\nu_n,\nu\ll dy$ then
    \[\mu_{n,x_{1:i}}\rightharpoonup \mu_{x_{1:i}}\Longrightarrow \mu_{n,x_{1:i}}^{i+1}\rightharpoonup\mu_{x_{1:i}}^{i+1}\]
and
    \[\nu_{n,y_{1:i}}\rightharpoonup \nu_{y_{1:i}}\Longrightarrow \nu_{n,y_{1:i}}^{i+1}\rightharpoonup\nu_{y_{1:i}}^{i+1}.\]
    By stability of optimal transport and uniqueness of narrow limits, $\gamma_{\KR,n}\rightharpoonup \gamma$, which is some plan that is the optimal plan from $\mu_{x_{i:1}}^{i+1}$ to $\nu_{y_{i:1}}^{i+1}$ with respect to cost $(x_{i+1}-y_{i+1})^2$. Doing this procedure for all $i\in\{2,\ldots,d-1\}$, we conclude that $\gamma$ is the optimal plan for the quadratic cost on each conditional measure of $\mu$ and $\nu$. This means that $\gamma=\gamma_\KR$, or that it is the Knothe-Rosenblatt plan from $\mu$ to $\nu$
\end{proof}

\section{Supporting results for Section~\ref{section3}}\label{appendix: C}

In the proofs of Section~\ref{section3}, we want to use the fact that $T_{\KR}$ is Lipschitz continuous. To be able to say that $T_{\KR}$ is Lipschitz, we must add more assumptions on our source and target measures. We also show some theorems on convergence of a sequence of functions.

\begin{lemma}[Brenier's theorem for rescaled quadratic cost, following \cite{baptista2024conditionalsimulationentropicoptimal}]\label{lem: brenier rescaled}
Suppose that $\mu\ll dx$ and $\nu$ has finite second moment. %
Then, the optimal map $\overline{T}$ pushing $\mu$ to $\nu$ according to cost $c(x,y) = \|A(x-y)\|^2$ for a positive definite matrix $A$ has the form 
\[
\overline{T}(x) = A^{-1} \Tilde{T}(Ax),
\]
where $\Tilde{T}$ is the optimal map according to the quadratic cost $\|x-y\|^2$ pushing measure $(A\cdot)_\sharp\mu=\mu(A^{-1}\cdot)$ to $(A\cdot)_\sharp\nu=\nu(A^{-1}\cdot)$.
\end{lemma}

\begin{theorem}[Caffarelli's regularity theorem]\label{thm: caf reg}
    Let $\Omega$ and $\Omega'$ be two bounded smooth open sets in $\R^n$. Let $\mu(dx)=f(x)dx$ and $\nu(dy)=g(y)dy$ be two probability measures such that $f = 0\in R^n\setminus \Omega$ and $g= 0\in \R^n\setminus \Omega'$. Assume that $f,g$ are both strictly positive and bounded away from infinity. If we further assume that $f,g$ are $C^\infty$ on their domains and $\Omega,\Omega'$ are uniformly convex, the optimal transport map (with respect to $\|x-y\|^2$) $T:\overline\Omega\to\overline{\Omega'}$ is a smooth diffeomorphism. 
\end{theorem}

Next, we present two corollaries of the above theorem.
\begin{corollary} \label{cor: TK Lip}
    Under the same conditions as in Theorem~\ref{thm: caf reg}:
    \begin{enumerate}
            \item The optimal transport map pushing $\mu$ to $\nu$ for the cost $c_\eps(x,y)$ is a smooth diffeomorphism;
            \item The Knothe-Rosenblatt map $T_{\KR}$ pushing $\mu$ to $\nu$ is Lipschitz continuous.
        \end{enumerate}
\end{corollary}

\begin{proof}
$ $
\begin{enumerate}
    \item By Lemma \ref{lem: brenier rescaled}, $T_\eps(x)= A_\eps \tilde T_\eps(A^{-1}_\eps x)$ where $A_\eps = \mathrm{diag}(1,\eps,\eps^2,\ldots,\eps^{d-1})$ and $\tilde T_\eps$ is the optimal map under $\|x-y\|^2$ pushing $(A_\eps)_\sharp\mu$ to  $(A_\eps)_\sharp\nu$. Since these two sets still satisfy the conditions on Theorem \ref{thm: caf reg} for fixed $\eps$, $T_\eps$ is a smooth diffeomorphism.
    \item We want to show that for any component $T_i(x_i|x_1,\ldots x_{i-1})$ of the $T_{\KR}$ map, the partial derivative $\frac{\partial }{\partial x_i}T_i(x_i|x_{1:i-1})$ is bounded. Recall that $T_i$ is the monotone rearrangement from $f(x_i|x_{1:i-1})dx_i$ to $g(y_i|y_{1:i-1})dy_i$. it follows that $G(T_i(x_i|x_{1:i-1})|T_{1:i-1}(x_{1:i-1})) = F(x_i|x_{1:i-1})$ where $F$ and $G$ the the conditional distribution functions of $f$ and $g$, respectively. By differentiating with respect to $x_i$, we have
    \[
    g(T_i(x_i|x_{1:i-1})|T_{1:i-1}(x_{1:i-1}))\frac{\partial }{\partial x_i}T_i(x_i|x_{1:i-1}) = f(x_i|x_{1:i-1}).
    \]
    It then follows that $\frac{\partial }{\partial x_i}T_i(x_i|x_{1:i-1})$ is bounded, since the densities are assumed to be bounded away from 0 and infinity. Since this is true for all components $i\in\{1,\ldots,d\}$, the Knothe-Rosenblatt map is Lipschitz.
\end{enumerate}    
\end{proof}

\begin{lemma}[Convergence of inverse maps adapted from Theorem 2 in \cite{barvinek1991convergence}]\label{lem: inverse conv}
Let $(M, d_M)$ and $(N, d_N)$ be metric spaces, with $N$ locally compact. 
Suppose $(f_n)$ is a sequence of bijective maps $f_n : M \to N$ converging almost uniformly to a function $f : M \to N$. Assume that $\mu$ is a finite measure on $M$, and that $\nu_n, \nu$ are finite measures on $N$ such that
\[
(f_n)_\sharp \mu = \nu_n, \quad \text{and} \quad f_\sharp \mu = \nu.
\] 
If each density of $\nu_n$ is bounded below by a positive constant $c(n) > 0$, and $f^{-1}$ is continuous, then
\[
f_n^{-1} \to f^{-1} \quad \text{$dy$-a.e. on } N.
\]
\end{lemma}
\begin{proof}
    Let $K_0$ be a compact set in $N$ and pick out another compact set $K_1$ such that $K_0\subseteq \mathrm{int}\ K_1 \subseteq K_1 \subseteq N$. Observe that since $f_n$ is almost uniform: For all $\eps>0$, there exists $E_\eps$ such that if $\mu(E_\eps)<\eps$ then $f_n \Rightarrow f$ on $E^c_\eps$. Set $S_{\eps,n} = f_n(E_\eps)$. Let $y\in K_0\setminus S_{\eps,n}$, this means that $y\notin f_n(E_\eps)$. It follows that $x_n:=f_n^{-1}(y)\in E_\eps^c$. Next, denote $y_n:= f(x_n)$ and choose $n$ sufficiently large such that $y_n\in K_1$. Now by the almost uniformness of $f_n$, for any $\eps>0$,
    \[
    d_N(f_n(x_n),f(x_n)) = d_N(y-y_n)\leq \eps
    \]
    for $n$ sufficiently large. Recall that since $f^{-1}$ is continuous, it is uniformly continuous on $K_1$. This means that for all $\eps>0$,
    \[
    d_M(f^{-1}(y)-f^{-1}(y_n))<\eps.
    \]
    Since $f^{-1}_n(y)=x_n = f^{-1}(y_n)$ since $y_n = f(x_n)$. It then follows that
    \[
    d_M(f^{-1}(y)-f^{-1}_n(y)) <\eps
    \]
    which is true for all $y\in K_0\setminus S_{\eps,n}$. Denoted $S_n := \bigcap_{k\in\N} S_{\frac1k,n}$. This set is in the $\sigma$-algebra so it is $\mu$-measurable. Notice by monotonicity, for each $k\in \N$
    \[
    \nu_n(S_n)\leq \nu_n(S_{\frac1k,n}) = \nu_n(f_n(E_{\frac1k})) = \mu(f^{-1}_n(f_n(E_{\frac1k}))) = \mu(E_{\frac1k})\leq\frac1k.
    \]
    Hence $\nu_n(S_n) = 0$. Consider some $y\in K_0$. Notice that if $y\in S_n^c$, then $y\notin S_{n,\frac{1}{k'}}$ for some $k'\in\N$ sufficiently large. We now want to get rid of the dependence on $n$ by performing a similar argument. Denote $S:= \bigcap_{n\in N} S_n$, which is measurable. It follows by monotonicity that and the fact that $\nu_n$
    \[
    c(n)dy(S)\leq \nu_n(S)\leq \nu_n(S_n) =0.
    \]
    Then if $y\in S^c$, then $y\notin S_{n'}$ for some $n'\in \N$ sufficiently large. Then clearly $y\notin S_{n',\frac{1}{k'}}$ for $n',k'\in\N$ sufficiently large. Thus $f^{-1}_n(y)\to f^{-1}(y)$ pointwise oustide of the null set $S$.
\end{proof}

\begin{lemma}\label{lem:xt opt}
Let $T_\eps$ be the optimal transport map corresponding to the weighted squared Euclidean cost $\|\cdot\|_\eps$ from $\rho_0=\mu$ to $\rho_1=\nu$ where $\mu\ll dx,\nu\ll dy$. By Brenier's theorem, $\gamma = (id\times T_\eps)_\sharp\mu$ is the optimal plan. Consider $\rho_t = (X_\eps(t))_\sharp\mu$ where $X_\eps(t,x)=(T_\eps(x)-x)t+x$. Then $X_\eps(t)$ is the optimal transport map from $\mu$ to $\rho_t$. 
\end{lemma}

\begin{proof}
     Let $\pi_t(x,y)=(y-x)t+x$, then 
    \begin{align*}
    (\pi_t)_\sharp\gamma &= (\pi_t)_\sharp(id\times T_\eps)_\sharp\mu\\
    &=(\pi_t\circ (id\times T_\eps))_\sharp\mu\\
    &=(X_\eps(t))_\sharp\mu.    
    \end{align*}
 Then let $\gamma_t := (\pi_0,\pi_t)_\sharp\gamma$, we get
\begin{align*}
    K_\eps(\rho_0,\rho_t)^{\frac12} &\leq \left(\int_{\R^d\times\R^d} \|z-z'\|_\eps^2 d\gamma_t(z,z')\right)^{\frac12}\\
    &=\left(\int_{\R^d\times\R^d} \|x-\pi_t(x,y)\|_\eps^2 d\gamma\right)^{\frac12}\\
    &=t\left(\int_{\R^d\times\R^d} \|x-y\|_\eps^2 d\gamma\right)^{\frac12} = tK_\eps(\rho_0,\rho_1)^\frac12.
\end{align*}
Applying this argument to the interval $[1-t,1]$, we deduce
\[
\begin{dcases}
    0\leq K_\eps(\rho_0,\rho_t)\leq tK_\eps(\rho_0,\rho_1)\\
    0\leq K_\eps(\rho_{1-t},\rho_{0})\leq (1-t)K_\eps(\rho_0,\rho_1).
\end{dcases}
\]
This allows us to say
\[
K_\eps(\rho_0,\rho_t)+ K_\eps(\rho_{1-t},\rho_{1})\leq K_\eps(\rho_0,\rho_1).
\]
The opposite inequality is obviously true, since $K_\eps(\rho_0,\rho_1)$ is the minimal cost of going from $\rho_0$ to $\rho_1$. Hence $K_\eps(\rho_0,\rho_t)+ K_\eps(\rho_{1-t},\rho_{1})= K_\eps(\rho_0,\rho_1),$ and all inequalities become equalities.
\[
 K_\eps(\rho_0,\rho_t)= \int_{\R^d\times\R^d} \|z-z'\|_\eps^2 d\gamma_t(z,z').
\]
Hence, $\gamma_t$ is optimal from $\rho_0$ to $\rho_t$, and $(id\times X(t))_\sharp\mu = ((\pi_0,\pi_t)\circ(id\times T))\sharp\mu = (\pi_0,\pi_t)_\sharp\gamma = \gamma_t$. Hence, $X_\eps(t)$ is the optimal transport map from $\rho_0$ to $\rho_t$.
\end{proof}

\bibliography{refs}{}

\bibliographystyle{alpha}

\end{document}